\documentclass[11pt,a4paper,leqno]{amsart}
\usepackage{amsmath,amsfonts}
\usepackage{amssymb}
\usepackage{upref}
\usepackage[mathcal]{eucal}
\usepackage{graphicx}

\date{\today}

\begin{document}
\title{The index and its prime divisors}

\author{Maciej P. Wojtkowski}
\address{ The University of Arizona \\
 Tucson, Arizona 85721, USA}
\email{maciejw@math.arizona.edu}

\date{\today}

            \theoremstyle{plain}
\newtheorem{lemma}{Lemma}
\newtheorem{proposition}[lemma]{Proposition}
\newtheorem{theorem}[lemma]{Theorem}
\newtheorem{corollary}[lemma]{Corollary}
\newtheorem{fact}[lemma]{Fact}
   \theoremstyle{definition}
\newtheorem{definition}{Definition}
\newtheorem{example}{Example}
    \theoremstyle{remark}
\newtheorem{remark}{Remark}

\newcommand{\macierz}[4]{\scriptsize{\left[\begin{array}{cc} #1 & #2 \\
        #3 & #4 \end{array}\right]}}

\maketitle

\begin{abstract}
  We propose a new interpretation of the classical
  index  of appearance for second order linear recursive
  sequences. It stems from the formula
  \[
   C_{n}(t)-2 =\frac{\Delta}{Q^{n}}\  L_n^2,\ \ \
   \text{where} \ \ t= (T^2-2Q)/Q, \  \Delta = T^2-4Q, 
   \]
   connecting the Chebyshev polynomials of the first kind
   $C_n(x)$ with the Lucas sequence defined for integer
   $T,Q\neq 0$ by the recursion $L_{n+1}= TL_n-QL_{n-1},
   L_0=0, L_1 = 1$. We build on the results of \cite{L-W}.

   We prove that for any prime $r\geq 2$ the sets 
   $\Pi_j(t,r), j=1,2,\dots$, of primes $p$
 such that $j$ is the highest power
 of $r$ dividing the index of appearance, have prime density
 equal to
 $\frac{1}{(r+1)r^{j-1}}$, for $r$-generic values of $t$.

 We give also complete enumeration of non-generic cases and
 the appropriate density formulas. It improves on the work
of  Lagarias, \cite{L}, and 
Ballot, \cite{B1},\cite{B2},\cite{B3}, on the sets of prime divisors
of sequences of ``finite order''.

Our methods are sufficient to prove that for any linear recursive
sequence of second order (with some trivial exceptions) the
set of primes not dividing any element
contains a subset of positive density.

We consider also some applications in arithmetic dynamics.

\end{abstract}

\section{Introduction}

Let us recall the setup from the paper \cite{L-W}.

For  rational  $ t$   let $D= D_{t} =
\left[ \begin{array}{cc} 0 &     -1 \\
    1  &   t \end{array} \right]$
be the matrix defining recursive sequences $\{x_n\}_{n\in\mathbb Z}$ by the
formula
\[
[x_n \  x_{n+1}]=[x_0 \  x_{1}]D^n.
\]
We choose to consider all rational initial conditions
$x_0,x_1$. It is traditional to study such sequences when they are integer
rather than
rational. Let us observe that our sequences are {\it almost integer} in
the sense that
all their elements can be represented as fractions with denominators
equal to $b a^n$ for fixed natural numbers $a$ and $b$. For such sequences
we study
the arithmetic properties of the sequence of numerators. The reason for
this approach
is not so much attempting greater generality but rather to achieve a
harmonious theory.

We will need an even more general environment of the parameter $t$ from
a field  $K$. In such a case we consider the ring
$\mathcal R = \mathcal R(t)=\mathcal R(t;K)$ of
$2\times 2$ matrices with
entries from $K$ which commute with $D=D_{t}$.
In particular $D \in \mathcal R$. A simple calculation
leads to the form of matrices $X \in \mathcal R$
\[X=
\left[ \begin{array}{cc} -x_{-1} &     -x_{0} \\
   x_{0} &   x_{1} \end{array} \right], \ x_1 = tx_0-x_{-1},
 x_0,x_1 \in K.
 \]
 For ease of notation
 we will denote an element $X \in \mathcal R$
by the second row $[x_0 \   x_1]$.
 Clearly if $ \delta = t^2-4$ is non-zero and  not a square in $K$ then
 the ring $K$ is isomorphic to the quadratic field extension
 $K(\sqrt{\delta})$.

In particular for a rational $t$ we will work with the ring 
$\mathcal R_p = \mathcal R_p(t)$ for any prime $p \geq 3$. It is obtained by
replacing $t$ with $t \mod p$, and all elements of the matrices
taken from the finite field $\mathbb F_p$. 
When using modular arithmetic with rational numbers we exclude
the finite set of prime divisors of the denominators of the fractions
involved.

Such rings in a more abstract fashion  were considered by Hall,
\cite{H}. Once we have the ring $\mathcal R$
we introduce the multiplicative subgroups $\mathcal S$
and $\mathcal S_p$  of matrices with determinant equal to $1$
with all elements from $K$ or $\mathbb F_p$,
respectively.
This group is reminiscent, but not identical to the Laxton group, \cite{Lx},
or the sequence group from \cite{L-W}. 
We will call it {\it the special sequence group}.

For the rest of the paper we assume that $t\neq 0,\pm 1, \pm 2$.
Moreover we define $\Pi(t)$ as the set of primes 
$p\geq 3$ not dividing the denominator of $ t$.
We have the following theorem.
\begin{theorem}
\label{Theorem1}
For $p \in \Pi(t)$ the group $\mathcal S_p$ is cyclic and
its order 
is equal to $p-1$ if $\delta \mod p$ is a square in $\mathbb F_p$,
and it is equal to $p+1$
if 
$\delta \mod p$ is not a square in $\mathbb F_p$.

If $t=2 \mod p$ then the order of $\mathcal S_p$
is equal to $p$ and if $t=-2 \mod p $ then it is equal
to $2p$. 
\end{theorem}
This theorem is well known in other formulations. 

\begin{definition}
 For $t\neq 0,\pm 1, \pm 2$ and $p\in\Pi(t)$
    {\it the index of appearance $\chi=\chi(t,p)$},
or {\it the index} for short, is the order of
the matrix $D_t \mod p$ in the special sequence group $\mathcal S_p$.
\end{definition}
Let us note that for a prime divisor $p$ of the numerator of $\delta = t^2-4$
the index is equal to $p$ or $2p$.

In other words $\chi$ is the smallest natural $k$ such that
$D^k = I \ \mod p$.
The powers $D^n$ for $n\in \mathbb Z$ can be expressed by  
the Chebyshev
polynomials of the second type $U_n(x)$, namely
\[
D^n = U_{n}(t)D-U_{n-1}(t)I =
\left[ \begin{array}{cc} -U_{n-1}(t) &     -U_{n}(t) \\
   U_{n}(t) &   U_{n+1}(t) \end{array} \right]=
\left[\  U_{n}(t) \ \   U_{n+1}(t)\ \right].
\]
Hence the index $\chi(t,p)$
is equal to the smallest natural number $n$ such that $U_{n}(t) = 0 \mod p$
and
$U_{n+1}(t)= - U_{n-1} = 1 \mod p$.

The trace  $ tr\ D^n = C_n(t)$, where $C_n(x) = U_{n+1}(x)-U_{n-1}(x)$
is the Chebyshev
polynomial of the first kind, $C_{0}(x) =2, C_1(x) =x$. For $p$ not dividing
the numerator of $\delta$ the index $\chi(t,p)$
can be defined equivalently as
the smallest natural number $n$ such that $C_{n}(t) = 2 \mod p$.

The choice of the name for the index $\chi$ is warranted
by the fact that it is equal to the classical index of appearance for
the {\it Lucas sequence} $\{L_n(T,Q)\}_{n\in \mathbb Z}$, \cite{E-P-Sh-W},
defined for
integer $T$ and $Q\neq 0$ by the recursive relation $L_{n+1} = TL_{n}-QL_{n-1}$
and the initial conditions $L_0=0, L_1=1$. The index of appearance $\mod p$
for the Lucas sequence is the smallest natural $k$ such that $L_k =0$.
It was established in \cite{W}(also \cite{L-W}, formula (1)),
that for $t= (T^2-2Q)/Q$  
\[
L_{2m}= TQ^{m-1}U_m(t), \ \ \ L_{2m+1}= Q^{m}W_{2m+1}(t), m = 0, \pm 1, \pm 2, \dots,
\]
where $W_{2m+1}=U_{m+1}+U_{m}$. We also introduce the polynomials
$V_{2m+1}=U_{m+1}-U_{m}$, to formulate a useful factorization
$U_{2m+1}=W_{2m+1}V_{2m+1}$, \cite{Y}. For even indices we have
$U_{2m}=C_{m}U_{m}$.
It follows from these formulas that for any integer $n$
\begin{equation}
  C_{n}(t)-2 =\frac{\Delta}{Q^{n}}\  L_n^2,\ \ \
  \text{where} \ \ \Delta = T^2-4Q.
\end{equation}
This formula can be also established directly by the  induction on $n$.
It is a consequence of (1) that the index of Definition 1
is equal to the classical index of appearance.


For a fixed prime $r$, including $r=2$, we remove $r$ from $\Pi(t)$
without changing the notation,
and introduce the partition
$\Pi(t) =  \bigcup_{j=0}^{+\infty}\Pi_j(t)$
where 
\[
\Pi_j(t)=\Pi_j(t,r) = \{p\ | \ r^j || \chi(t,p)\}.
\]
Let us note that the prime divisors of the numerator of $\delta=t^2-4$
belong to $\Pi_0$ for $r\geq 3$, and to $\Pi_0\cup \Pi_1$ for $r=2$.

Our main result is the following theorem.
Depending on the prime $r$ we will 
additionally exclude certain values of $t$, to be described later.
The remaining values
of $t$ are called {\it r-generic}. 
\begin{theorem}
\label{Theorem2}
For any prime $r\geq 2$ and an $r$-generic $t$ the disjoint subsets of primes
$\Pi_j(t,r), j= 1,2,\dots,$
have the prime densities and they are equal to
\[
|\Pi_j| = \frac{1}{(r+1)r^{j-1}}.
\]
\end{theorem}
Clearly also $\Pi_0$ has the prime density $|\Pi_0|= 1-\frac{r}{r^2-1}$.
The assumption of genericity includes the condition that $t^2 -4$
is not a rational square. In Section 11 we will outline the necessary
changes in formulations and proofs required when $t^2-4$ is a rational
square, the case called {\it reducible}.

This theorem covers the calculation of prime densities of the set of prime
divisors of
most recursive sequences of ``finite order''
studied in the past, and that only
for $r=2$ and $r=3$. This problem was first addressed
in ground-breaking papers of Hasse, \cite{H} and Lagarias, \cite{L}.
It was followed by the extensive work of 
Ballot, \cite{B1},\cite{B2},\cite{B3}. More results can be
found in papers of Moree,\cite{M1},\cite{M2},  Moree and Stevenhagen,\cite{M-S},
and the unpublished thesis of Lujic,\cite{Lj}.
The survey paper of Moree \cite{M3} gives a summary of these advances
in Chapter 8.4.

After proving Theorem 2 in Sections 2,3 and 4
we proceed with complete investigation of the non-generic
cases which was not accomplished before. The
intensively studied case of $r=2$ is the most involved.
The first condition of
genericity is primitivity, namely $t$ is $r$-primitive if
there are no rational solutions to the equation $C_r(x)=t$.
In Proposition 11 we show the simple connection of the
partitions for $t$ and $t_r=C_r(t)$.

For $r\geq 3$ the non-genericity appears when $t^2-4 = \pm r b^2$
for a rational
$b$,  with $+$ for $r=1 \mod 4$ and $-$ for $r=3 \mod 4$.
In the  $+$ case a slight modification of the proof delivers
the formula from Theorem 2. In the $-$ case we have Theorem 13
where the densities are doubled compared with Theorem 2.

Further we formulate
three types of symmetry for the problem, discovered
in \cite{L-W}. The first one
is the universal {\it twin} symmetry $t\leftrightarrow -t$. The
{\it circular} symmetry occurs when $t^2-4=-b^2$ for a rational $b$,
and the {\it cubic} symmetry when  $t^2-4=-3f^2$ for a rational $f$.
In these two cases $\mathcal R(t)$ is isomorphic 
to the cyclotomic fields of order $4$ and of order $3$, respectively. Special
finite order sequences 
in these cases were studied by Ballot, \cite{B3}.
In all three cases there are {\it associate values} for a given
$t$, which are the solutions of the equations $C_k(x)=C_k(t)$,
for $k=2,4$ and $3$, respectively. This phenomenon has an impact on
the definition of $r$-genericity for $r=2$ (twin and circular symmetry)
and $r=3$ (cubic symmetry). In these cases  the genericity conditions
includes the
requirement that also associate values are $2$-primitive for $r=2$
and $3$-primitive for $r=3$. These 
conditions are called twin primitivity, circular primitivity
and cubic primitivity, respectively.
For every $t'$ there is always $t$ satisfying this more stringent
primitivity condition such that $C_{r^m}(t)$, for some natural $m$,
is equal to $t'$ or to an  associate value of $t'$.
Proposition 11
is then augmented in Corollary 24, (7) and (8), and Proposition 14 and
(6), respectively, giving the connection between the partitions for
$t$ and $t'$.  It remains to give the densities for the cubic case
(Theorem 13 and Section 6), the circular case (Theorem 28) and two cases
of non-genericity for $r=2$: {\it type A}, when $2(2\pm t)$
is a rational square
(Theorem 27), and {\it type B}, when $2(t^2-4)$ is a rational square
(Theorem 25).

We employ a  common inductive scheme 
for  all the density 
calculations in terms of the special sequence groups
$\mod p$. It is  introduced in Section 2. In Section 3
we establish for $r\geq 3$  an equivalent description in terms of
the factorization of appropriate polynomials $\mod p$,
setting the stage for the application of Frobenius-Chebotarev
density theorem, \cite{S-L}. The case of $r=2$ is treated in Sections
7, 8 and 9. It starts
with the introduction of modified polynomials (Theorem 17) in the
characterization of divisibility properties by powers of $2$.
The rest of the effort involves calculations of the orders
of respective Galois groups, which are different in generic
and non-generic cases.






In the cubic case there are always sequences of
order $3$. The densities for their set of divisors
were obtained  by
Lagarias, \cite{L},  and Ballot, \cite{B3}.
They are covered by Theorem 13 for $r=3$.
In the same scenario there are also sequences 
of order $6$ with the set of divisors
characterized by the divisibility of the index by $6$.
Ballot established
their prime density, \cite{B3}, under some genericity condition.
Using our approach we obtain some improvement of this result in Section
10. We show in Theorem 33 that under the conditions
of cubic primitivity and twin primitivity the divisibility 
properties of the index by powers of $2$ and by powers of $3$
are ``independent''.

Ballot, \cite{B1}, introduced an interesting sequence, 
for any $r\geq 2$,
\[
B_k(T,Q) := \frac{L_{rk}(T,Q)}{L_{k}(T,Q)}, \ k = 1,2,\dots,
  \]
  where $\{L_n(T,Q)\}_{n\in \mathbb Z}$ is the Lucas sequence.
This sequence is second order linear only for $r=2$.
Ballot proved that the set of prime divisors
of the sequence is equal to $\Pi(t) \setminus \Pi_0(t,r)$.
In Section 10 we give another proof 
using the above expression of the Lucas sequence by the Chebyshev polynomials.
  By Theorem 2 the prime
  density of the set of prime divisors is 
  equal to $\frac{r}{r^2-1}$ (in the $r$-generic case).
  Ballot established such a formula
  in \cite{B1} in the reducible case only, when
  $\Delta=T^2-4Q$ is a rational square.

For any $r\geq 2$ there are
  linear recursive sequences of second order  with the set
  of prime divisors equal to $\Pi_0(t)$ for $r\geq 3$ and
  $\Pi_0(t)\cup\Pi_1(t)$ for $r=2$. We will discuss them
  in Section 10.

  In Section 11 we give three applications of Theorem 2,
  and its non-generic versions, in arithmetic dynamics.
  The first one is about rational orbits of irrational rotations.
  These orbits are given by the linear recursive sequences.

  The second application is the study of prime divisors of rational orbits
  of  chaotic interval maps given by the
  Chebyshev polynomials $C_k(x):[-2,2] \to [-2,2]$. The claim is
  that for any rational initial condition the set of prime divisors
  of the resulting orbit (which is an almost integer sequence)
  has density zero. Rafe Jones, \cite{J}, proved such a claim for
  a large family of quadratic functions, including $C_2(x) =x^2-2$.
  However, he restricted the study to integer initial values.
  We conjecture that his result is valid also for rational initial
  conditions, and the resulting almost integer orbits.

  In the  last application we make a connection between Theorem 2 and
  the study of iterations of the quadratic polynomial
  $f_t= (x+t)^2-2-t$ by Gottesman and Tang,\cite{G-T}.
  They found the prime density of the set of prime divisors
  of the orbit of $x=0$ and in contrast to the result of Jones,
  it is positive. They restricted their attention to integer
  values of $t$. We show that their problem is covered 
  by Theorem 2, and $t$ can be any rational value with appropriate
  exclusions.

  In Section 12 we employ the methods
  developed in the proof of Theorem 2
  to prove  for any linear recursive sequence
  of second order, with the exception
  of the identity, i.e, the Lucas sequence,   
  that the set of primes not dividing any element
  contains a subset of positive density.    This is
  a strengthening of the theorem of Schinzel, \cite{S1},
  who established in the reducible case the existence of infinitely
  many non-divisors.

  Although there is no direct connection with the present paper
  the author derived original inspirations from the papers of
  Schur, \cite{Sch}, Rankin, \cite{R}, and Arnold, \cite{A}.

  \section{Proofs of Theorems 1 and 2 (part one)}
  \begin{proof}{(of Theorem 1)}
    If $\delta \mod p $ is a square in $\mathbb F_p$ then
    the matrix $D \mod p$ has two different eigenvalues in
    $\mathbb F_p$. Hence it can be diagonalized over $\mathbb F_p$
    by a conjugation. The same conjugation diagonalizes also all the matrices
    that commute with $D$. Hence we can naturally identify $\mathcal S_p$ with
    diagonal matrices with determinant $1$. Clearly there are $p-1$
    such matrices and they form a group naturally isomorphic to
    $\mathbb F_p^*$, which is cyclic.

    If $\delta \mod p $ is not a square in $\mathbb F_p$ then
    the ring of matrices commuting with $D \mod p$ is a finite
    field, isomorphic to $\mathbb F_{p^2}$. It is easy to see that
    the determinant of a matrix corresponds to the norm in  
    $\mathbb F_{p^2}$ over $\mathbb F_{p}$. Hence the group $\mathcal S_p$ is
    isomorphic to the subgroup of elements of $\mathbb F_{p^2}$
    of norm $1$, i.e., the solutions  of the  equation $x^{p+1} =1$
    in $\mathbb F_{p^2}$. Since $\mathbb F_{p^2}^*$ is cyclic of order
    $p^2-1$ then it contains $p+1$ different roots of unity of degree $p+1$.

    It can be
    readily established that the recursive sequence 
    with the initial
    conditions $x_0 = 0, x_1= 1$  is equal to $x_n = n$ for $n\geq 0$ when
    $t= 2 $, and to $x_n = (-1)^{n+1}n$ for $n\geq 0$ when
    $t= -2 $.
    It follows that if $t^2=4 \mod p$ then
    $D^p = [ x_p \ x_{p+1} ] = \pm I \mod p$.
    Hence the order of $D$ in $\mathcal S_p$ is equal to $p$ or $2p$.
    \end{proof}

  To prove Theorem 2 we begin by developing an inductive
  description of the elements of the partition.  
To ease the notation, for a given $t$ and an odd prime $p$,
we denote by  $\widehat p= \widehat p(t,p)$ the order of the
group $\mathcal S_p$.
By Theorem 1 $\widehat p = p\pm 1$, except for the divisors of $\delta$.
In the following we will exclude from consideration these exceptional divisors:
they always  belong to $\Pi_0$ (or  $\Pi_0 \cup \Pi_1$ for $r=2$),
and anyway since there are only finitely many of them they cannot affect
the densities.
With a fixed  prime $r$, including
$r=2$, let
\[
\mathcal K_j = \{ p \in \Pi(t)\ | \ r^j |\widehat p\}, \ \
\mathcal M_j = \{ p \in \Pi(t)\ | \ \exists Y\in \mathcal S_p(t),  Y^{r^j}=D\},
\]
for $j=1,2, \dots$, and 
$\mathcal M_0 =\mathcal K_0= \Pi(t)$. 
Let us  note that if $p \notin \mathcal K_1$ then $p\in \mathcal M_n$
for any natural $n$. It does not play any role in the following
because we consider only the subsets $\mathcal M_n\cap \mathcal K_j$
with $j\geq n$. Clearly for $r=2$ we have $\mathcal K_1 = \Pi(t)$
for any $t$. We will later find out that the same conclusion holds also
for $r=3$ and cubic $t$.

Since the order of any element of a finite group divides the order 
of the group we get that if $p \notin \mathcal K_1$ then $p \in \Pi_0(t,r)$.


We split $D_0 =: \mathcal K_1$ into four disjoint subsets,
\[
\begin{aligned}
  A_1 &=\left(\mathcal M_0\setminus \mathcal M_1\right)\cap
  \left(\mathcal K_1\setminus \mathcal K_2\right),\ 
B_1=\mathcal M_1\cap \left(\mathcal K_1\setminus \mathcal K_2\right),\\
C_1&=\left( \mathcal M_0\setminus \mathcal M_1\right)\cap \mathcal K_2,\
\ \ \ \ \ \ \ \ \ \ D_1=\mathcal M_1\cap \mathcal K_2.
\end{aligned}
\]
\begin{proposition}
  The following inclusions hold
  \[
  B_1\subset \Pi_0, \ \ A_1 \subset \Pi_1,\ \
  C_1\cap\left(\mathcal K_j\setminus \mathcal K_{j+1}\right) \subset \Pi_j,
\]
for $j= 2,3,\dots$. 
  \end{proposition}
This splitting and the Proposition 3 are illustrated   
in the following Table 1.

\begin{table}[htbp]
  \caption{The splitting of $\mathcal K_1$ and
    $D_1=\mathcal M_1\cap\mathcal K_2$ into four disjoint subsets }  
  $$
  \begin{array}{|c|c|c|}
\hline
& \mathcal M_0\setminus \mathcal M_1
& \mathcal M_1\\
 \hline
 \mathcal K_1\setminus \mathcal K_2  & r \parallel \chi & r \nmid \chi \\
\hline
\mathcal K_2   &   r^2 \mid \chi                 &  ? \\
 \hline
  \end{array}\ \ \
  \begin{array}{|c|c|c|}
\hline
&\mathcal M_1\setminus \mathcal M_2
& \mathcal M_2 \\
 \hline
\mathcal K_2\setminus \mathcal K_3  & r \parallel \chi & r \nmid \chi  \\
\hline
\mathcal K_3   &    r^2 \mid \chi  &  ? \\
 \hline
  \end{array}
  $$
\end{table}
\begin{proof}
  By Theorem 1 the group $\mathcal S_p(t)$
  has a generator $Z$. Without loss of generality we can assume
  that $D= Z^s$ with $\chi = \widehat p/s$.
  If $p \in B_1$ then $r \mid s$ and $r^2 \nmid \widehat p$. It follows that
  $p \in \Pi_0$.

  Similarly, if
  $p\in A_1\cup C_1= \mathcal M_0\setminus \mathcal M_1 $ then $r \nmid s$.
  Hence if $p\in A_1$ then
  $r \parallel \widehat p$ and
  $p \in \Pi_1$. If
  $p\in C_1\cap \left(\mathcal K_j\setminus\mathcal K_{j+1}\right)$
  then 
  $r^j \parallel \widehat p$ and
  $p \in \Pi_j$, for $j\geq 2$ . 
\end{proof}
We proceed defining inductively the splitting of the set $D_{n-1}=
\mathcal M_{n-1}\cap\mathcal K_{n}$
into four different subsets $A_{n},B_{n},C_{n},D_{n}$.
\[
\begin{aligned}
  A_n &=\left(\mathcal M_{n-1}\setminus \mathcal M_n\right)\cap
  \left(\mathcal K_n\setminus \mathcal K_{n+1}\right),\ 
B_n=\mathcal M_n\cap \left(\mathcal K_n\setminus \mathcal K_{n+1}\right),\\
  C_n&=\left( \mathcal M_{n-1}\setminus \mathcal M_n\right)\cap
  \mathcal K_{n+1},\
\ \ \ \ \ \ \ \ \ \ D_n=\mathcal M_n\cap \mathcal K_{n+1}.
\end{aligned}
\]
We get immediately the extension of Proposition 3
\begin{proposition}
  For any $n\geq 1$ and $m\geq 2$ the following inclusions hold
  \[
  B_n\subset \Pi_0, \ \ A_n\subset \Pi_1,\ \
  C_n\cap\left(\mathcal K_{n+m-1}\setminus \mathcal K_{n+m}\right)
  \subset \Pi_{m}.
  \]
  \end{proposition}
The proof does not differ substantially from the proof of Proposition 3.

Let us observe that to obtain the prime densities of all the sets
$A_{n},B_{n},C_{n},D_{n},$  $n=1,2,\dots$ it is enough to establish them
for the sets $\mathcal M_{n}\cap \mathcal K_{j}$, for all $n\geq 0$
and only for $j=n,n+1,n+2$. However, to get also the densities of the sets
$C_n\cap \Pi_m, m= 2,3,\dots$ we need to cover all values of $j\geq n$.

\section{Description of the sets $\mathcal M_n\cap \mathcal K_{j}$}

We are going to describe the sets of primes
$M_n\cap \mathcal K_{j}$
in terms of some polynomials
factoring over finite fields, thus allowing the application of the
Frobenius-Chebotarev Density Theorem.

For $j=1,2,\dots$ let us denote by $\Phi_j$ the $r^j$-th
cyclotomic polynomial 
\[
\Phi_j(x) = \frac{x^{r^j}-1}{x^{r^{j-1}}-1}.
\]

Further we put for a rational $t$ and  $n=1,2,\dots$,
\[
f(x) = x^2-tx+1, \ \  g_n(x) = C_{r^n}(x)-t.
\]
The polynomials $g_n$ are called   {\it DeMoivre's
  polynomials}, \cite{G}. We will say that a polynomial {\it splits
linearly} if it factors into linear factors, and that it {\it splits
   quadratically} if it factors into irreducible quadratic factors.
At this point we need to exclude $r=2$.
\begin{theorem}
  For any prime $r\geq 3$ and a natural $j$
  the set $\mathcal K_j= \mathcal K_j(t,r)$
  contains exactly the primes
  $p\in \Pi(t)$  such that
  either the polynomial $f\Phi_j\mod p$ splits linearly or it splits
  quadratically.
\end{theorem}
We will be using repeatedly the  following lemma characterising
the traces and determinants of the elements of $\mathcal R(t, K)$.
\begin{lemma}
  For any $a,b \in K, \ a^2-4b \neq 0$
  there is a matrix $Y \in \mathcal R(t,K)$
  such that $tr \ Y = a$ and $ \det Y = b$ if and only if
  $(t^2-4)(a^2-4b)$  is a square in the field of scalars $K$.
  Moreover for a given pair $a,b$ there are precisely two such matrices
  $Y$ and $(\det Y) Y^{-1}$.
\end{lemma}
\begin{proof}
  For $Y = [y_0\ y_1] \in \mathcal R(t,K)$ we get
  \[
  tr \ Y = 2y_1-ty_0 = a, \ \ \det Y = y_1^2-ty_1y_0 +y_0^2 = b.
  \]
  Substituting the first equation into the second we obtain
\[
(t^2-4)y_0^2 = a^2-4b.
\]
    \end{proof}
\begin{proof}{(of Theorem 5)}
  Clearly $f\mod p$ splits if and only if $t^2-4\mod p$ is a square in
  $\mathbb F_p$. If $f\mod p$ splits then
  the  matrix $D\mod p$ has two eigenvalues and two eigenvectors, so that 
  the ring $\mathcal R_p(t)$ is canonically isomorphic
  to the ring of diagonal matrices. Hence the special sequence group
  $\mathcal S_p(t)$ contains a primitive $r^j$-th root of unity if and only if
  $\mathbb F_p$ contains a primitive $r^j$-th root of unity. But that is
  equivalent to $\Phi_j\mod p$ having a root. It is well known
  that if $\Phi_j$ has a root then it splits linearly.

  If $f\mod p$ does not split then the ring 
  $\mathcal R_p(t)$ is canonically isomorphic to the finite
  field $\mathbb F_{p^2}$. If $r^j \mid \widehat p$ then $\mathcal S_p(t)$
  contains a matrix $A$ which is  a primitive root of unity of degree $r^j$.
  Hence $A$ is a root of $\Phi_j  \mod p$, and so $\Phi_j$  splits
  linearly over the field $\mathcal R_p$. The roots of $\Phi_j$
  are powers of $A$. For the roots $A^s$ and $A^{-s}$ we have
  \[
    \left(x-A^s\right)\left(x-A^{-s}\right) = x^2-a_sx+1, \ \
    a_s = tr \ A^s,
    \]
    where we tacitly identify the scalar multiplicities of $I$ in $\mathcal R_p$
    with elements of $\mathbb F_p$.
    By Lemma 5 $a_s^2-4$ is not a square in $\mathbb F_p$ and hence
    $x^2-a_sx+1$ does not split over $\mathbb F_p$.
    It follows that $\Phi_j$ splits
    over $\mathbb F_p$ into irreducible quadratic factors.

    Conversely, if $f$ is irreducible $\mod p$ and  $h = x^2-ax+b$
    is an  irreducible factor of $\Phi_j\mod p$,
    then $(t^2-4)(a^2-4b)$ is a square $\mod p$ and
    by Lemma 5 there is a matrix $Y\in\mathcal R_p$ such that
    $tr \ Y = a, \det Y = b$. Since $Y$ is a zero of $h$ it is also
    a zero of $\Phi_j$, so that $Y^{r^j}=I$. It follows that $r^j$ divides
    the order of the multiplicative group of $\mathcal R_p(t)$, equal to
    $p^2-1$. Since $r\geq 3$ we conclude that either $r^j\mid p-1$ or
    $r^j\mid p+1$. The former would imply that $\mathbb F_p$ contains
    a root of unity of degree $r^j$, and consequently linear splitting
    of $\Phi_j\mod p$. This contradiction implies the latter, which means
    that $p \in \mathcal K_j$.
\end{proof}
\begin{remark}In the last proof it was not necessary to argue that the
matrix $Y$ must belong to $\mathcal S_p$, i.e., by necessity $b=1$.
Let us do it for the sake
of clarity. We have established that $p+1 = kr^j$ for some
    natural $k$.   
    The other zero of $h$ is $\widetilde Y = Y^p$ and we get
    $bI=Y\widetilde Y = Y^{p+1}$. Since
    $Y^{p+1} = \left(Y^{r^j}\right)^k = I$ we conclude that $b=1$,
    i.e., $Y\in S_p(t)$.     
\end{remark}

We proceed with the characterization of the sets
$\mathcal M_n\cap \mathcal K_{j}$ for $j\geq n$. So far we have treated
the case $n=0$ and all values of $j\geq 0$.
\begin{theorem}
  For a prime $r\geq 3$ and any natural $n$ and $j\geq n$
  the set $\mathcal M_n\cap \mathcal K_{j}$
  contains exactly the primes
  $p\in \Pi(t)$  such that the polynomial $g_n\mod p$ splits linearly and
  either the polynomial $f\Phi_{j}\mod p$ splits linearly or it splits
  quadratically.
  \end{theorem}
We need another lemma  about 
the traces of elements in $\mathcal S$ or $\mathcal S_p$.
\begin{lemma}
  If $a\in K, \ a\neq \pm 2$, is a trace of a matrix in $\mathcal S$ and
  $\widetilde a$ is a zero of $C_m(x)-a$ in the field $K$,
  for a natural $m$, then there are matrices  in $\mathcal S$
  with the traces $\widetilde a$
  and $C_s(a)$, for any natural $s$.
\end{lemma}
The subject matter of this lemma is reminiscent of the work of
Rankin,\cite{R}, and in fact when combined with Theorem 1 it provides
the proof of the main result of the paper.
\begin{proof}
  We recall another identity
  for the Chebyshev polynomials,\cite{W}, for any natural $s$
  \begin{equation}
C_s^2(x) - (x^2-4)U_s^2(x) = 4.
\end{equation}
It implies that both $(C_s^2(a)-4)(a^2-4)$ and $(\widetilde a^2-4)(a^2-4)$
are squares  in $K$. The application of Lemma 5 finishes the proof.
    \end{proof}
\begin{proof}{(of Theorem 7)}
  If $p \in \mathcal M_n$ then there is a matrix $Y\in \mathcal S_p$ such
  that $Y^{r^n}=D$. Hence the trace $tr\ Y^{r^n} = t$. It follows that
  $y = tr\ Y$ is a zero of the DeMoivre's polynomial $g_n(x) = C_{r^n}(x)-t$.

  If $p \in \mathcal K_n$ then there is a matrix $Z\in \mathcal S_p$
  which is a primitive root of unity of order $r^n$. It follows that
  \[
  \left(YZ^k\right)^{r^n} = D, \ \ \text{for} \ \ k = 0,1,\dots,r^{n}-1. 
  \]
  Hence the traces $y_k= tr\ YZ^k$
  are different zeroes of the same DeMoivre's polynomial $g_n$.
  Indeed, should $y_k = y_l, k\neq l$ then by Lemma 5 either 
  $YZ^k = YZ^l$ or  $YZ^k = \left(YZ^l\right)^{-1}$. The first equality is
  impossible because $Z$ is a primitive root of unity. Raising both
  sides of the second equality to the power $r^n$ we get $D^2=I$.
  This implies that $2= tr\ D^2 = t^2-2$ and  hence $t^2 = 4 \mod p$,
  and we have excluded such values of $p$ from considerations.

  Thus we  arrive at the conclusion that the DeMoivre's
  polynomial $g_n\mod p$ splits linearly.  Together with Theorem 5
  this delivers the forward claim of the theorem.

  Conversely, if the polynomial $g_n \mod p$ splits linearly
  then
  by Lemma 7  there is a matrix $Y\in \mathcal S_p$ such that 
  $tr\ Y = y$, where $y \in \mathbb F_p$ is a zero of $g_n \mod p$.
  Further $tr\ Y^{r^n} =C_{r^n}(y) =t$ and
  hence $Y^{r^n}$ is equal to  $D$ or $D^{-1}$. Changing $Y$ to $Y^{-1}$,
  if needed, we can claim that $Y^{r^n}$ is equal to  $D$.
  Hence $p\in \mathcal M_n$. If additionally $f\Phi_{j}$ splits
  linearly or quadratically then by Theorem 5 $p \in \mathcal K_{j}$.
\end{proof}

\section{The Galois groups of $f\Phi_{j}g_n$}

We are going to study the splitting fields $H_{n,j}$ of the polynomials
$f\Phi_{j}g_n$, for all $n\geq 1, j\geq n$.
We are going to develop a different
interpretation of these fields. Let $\xi,\xi^{-1}\in\mathbb C$ be the
zeroes of $f$, and let 
$z_n\in\mathbb C$ be a root of $\xi$ of order $r^n$. Further for
any natural $k$ let 
$\zeta_k$ be a primitive root of unity in $\mathbb C$ of order $r^k$

  We have the following tower of field extensions
  \[
  \mathbb Q < \mathbb Q\left(\zeta_{j}\right)
  < \mathbb Q\left(\xi,\zeta_{j}\right) <
  \mathbb Q\left(\xi,\zeta_{j},z_{n}\right).
  \]
  \begin{theorem}
    For any rational $t \neq 0,\pm 1, \pm 2$ and any $n\geq 1$ and $j\geq n$
    \[
H_{n,j}=\mathbb Q\left(\xi,\zeta_{j},z_{n}\right).
\]
  \end{theorem}
  \begin{proof}
We will first prove the inclusion 
$H_{n,j} \subset \mathbb Q\left(\xi,\zeta_{j},z_{n}\right)$.
We have that $\zeta_{j}^{r^{j-n}}$ is a primitive root of unity of degree $r^n$,
so that we may as well denote it by $\zeta_{n}\in H_{n,j}$. Further
$z_n\zeta_{n}^{k},  k = 0,1,\dots, r^{n}-1$ are all the different zeroes of
the polynomial $x^{r^n}-\xi$. It follows from the identity
$z^s+z^{-s} = C_s(z+z^{-1})$ that
\[
y_k = z_n\zeta_{n}^{k}+z_n^{-1}\zeta_{n}^{-k},  k = 0,1,\dots,
r^{n}-1
\]
are all the different zeroes of $g_n(x) = C_{r^n}(x)-t$. Indeed
\[
C_{r^n}(y_k) = \left(z_n\zeta_{n}^{k}\right)^{r^n}
+\left(z_n^{-1}\zeta_{n}^{-k}\right)^{r^n} = \xi +\xi^{-1} = t.
\]
The inclusion is established.

Conversely, if $y_0$ is a zero of $g_n(x)$ then a zero of
$x^2-y_0x+1$ is a zero of $x^{r^n}-\xi$. The solutions of the
quadratic equation belong to $H_{n,j}$. Indeed by (2) we have
\[
\left(y_0^2-4\right)U_{r^n}^2(y_0) = C_{r^n}^2(y_0) - 4 = t^2-4 = (2\xi+t)^2.
\]
We thus conclude that
$\mathbb Q\left(\xi,\zeta_{j},z_{n}\right)\subset H_{n,j}$.
    \end{proof}

  Let us note that so far we have not restricted the  values of $t$
beyond $t\neq 0,\pm 1, \pm 2$. Many more restrictions will appear in
this Section with the introduction of $r$-genericity. At this stage
we assume additionally that $r \geq 3$. The case $r=2$ is much harder,
and it covers the density of the prime divisors of
the companion Lucas sequence.
The partition of primes for $r=2$ was studied in \cite{W}.
We will give in Sections 7,8,9 a new rendition of the main results
of that paper
using the simplifications introduced in \cite{L-W} and the current paper.
  \begin{definition}
    A rational $t\neq 0,\pm 1, \pm 2$ is $r$-primitive, for $r\geq 2$,
    if $C_r(x) -t$ has no rational roots.

    A rational $t\neq 0,\pm 1, \pm 2$ is $r$-generic, for $r\geq 3$,
    if it is $r$-primitive and
    $\xi\notin \mathbb Q(\zeta_1)$,
    i.e., $\sqrt{t^2-4}$ is not an element of the  cyclotomic field
    of order $r$.
  \end{definition}
  It is well known that the last condition is equivalent to
  $t^2-4 \neq \pm ra^2$ for a rational $a$, where $+$ is used for
  $r = 1\mod 4$ and $-$ for   $r = 3\mod 4$.
  This condition appears explicitly in the study by Girstmair,\cite{G},
  of the Galois groups for the DeMoivre's polynomials. 
  \begin{theorem}
    If a rational $t$ is $r$-generic, for $r\geq 3$ then for $j\geq n\geq 1$
    \[
    \left[\mathbb Q(\zeta_j):\mathbb Q \right] = (r-1)r^{j-1},
    \left[\mathbb Q(\xi, \zeta_j):\mathbb Q(\zeta_j) \right] = 2,
    \left[\mathbb Q(\xi, \zeta_j, z_n):\mathbb Q(\xi,\zeta_j) \right] =
    r^n.
    \]
    \end{theorem}
  \begin{proof}
    The first equality is just the degree of the cyclotomic extension of
    order $r^j$. The second follows from the assumption of $r$-genericity,
    that $\xi$
    does not belong to the cyclotomic extension of order $r$. Indeed
    the cyclotomic extension of order $r^j$ contains no more quadratic
    irrationals, beyond $\sqrt{\pm r}$.

    To prove the third equality let us observe that for $r$-primitive $t$
    the binomial $x^{r^n}-\xi$ has no roots in $\mathbb Q(\xi)$, or
    equivalently it is irreducible over $\mathbb Q(\xi)$.
    Indeed the ring $\mathcal R(t)$ is isomorphic to
    $\mathbb Q(\xi)$ and the isomorphism takes the matrix $D$ to $\xi$.
    As we have argued in the proof of Theorem 7,
    $g_1(x) = C_{r}(x)-t$ has a rational root if and only if
    $x^{r}-D$ has a root in $\mathcal R(t)$. Hence $r$-primitivity of $t$
    implies 
    irreducibility of  $x^{r}-\xi$ over $\mathbb Q(\xi)$, and clearly
    the same holds for all $x^{r^n}-\xi, n\geq 1$.

    The Galois group
    $Gal\left(\mathbb Q(\xi,\zeta_j)/\mathbb Q(\xi)\right)$
    is abelian. Moreover $\mathbb Q(\xi)$ does not contain
    an $r$-th root of unity. It follows from the fact that only
    the cyclotomic field of order $3$ is quadratic, and we explicitly
    ruled out in $r$-genericity that $\xi$ belongs to that field
    in case of $r=3$.
    Hence by Kummer theory, since  $x^{r^n}-\xi$ is irreducible
    over $\mathbb Q(\xi)$ then it is irreducible over its abelian field
    extension
    $\mathbb Q(\xi,\zeta_j)$. We conclude that  the radical extension
    in the third equality is indeed of degree $r^n$. 
    \end{proof}

  We are ready to finish the proof of Theorem 2 for $r\geq 3$.
  The case of $r=2$ will be treated in Section 7.
  \begin{proof}(the completion of the proof of Theorem 2 for $r\geq 3$)

    Theorem 10 allows an immediate application of the Frobenius-Chebotarev
    Density Theorem, \cite{S-L},
    to the calculation of the density of the set of primes
    for which
    the polynomial $g_nf\Phi_{j}, j \geq n$, splits linearly. We get that
    the prime density of this set is equal to
    $1/\left(2(r-1)r^{n+j-1}\right)$.

    It remains to address the same question in the case of the quadratic
    splitting of $f\Phi_j$. To that end we need more detailed information
    about the Galois group $Gal\left(H_{n,j}/\mathbb Q\right)$.
    We will use the interpretation of this group as a group of 
    permutations of the zeroes of $g_nf\Phi_{j}$.
      We are going to prove that 
       there is exactly one automorphism of $H_{n,j}$ which
      splits the zeroes of $f\Phi_j$ into $2$-cycles, and fixes
      the zeroes of $g_n$.

    Let us consider 
    an automorphism $\gamma$ of the cyclotomic field $\mathbb Q(\zeta_j)$.
    We have
    $\gamma(\zeta_j) = \zeta_j^s$, for some natural
    $s, 1 \leq s < r^j, r\nmid s $. If $\gamma$, as a permutation,  
    splits all zeroes of $\Phi_j$ into $2$-cycles, then  
$\zeta_j=\gamma(\zeta_j^s) = \zeta_j^{s^2}$. It follows
    that $\zeta_j^{s^2-1}=1$, so that $r^j\mid s^2-1$, and hence,
    since $r\geq 3$, either $r^j \mid s-1$ or $r^j \mid s+1$. The former
    is impossible and the latter leads to $s =r^j-1$,
    i.e., $\gamma(\zeta_j) = \zeta_j^{-1}$. Hence this is our only choice
    for the automorphism $\gamma$ of the cyclotomic field.

    The Galois group of the
    splitting field of $f\Phi_j$ is the direct product of the
    Galois groups of $f$ and of $\Phi_j$. Hence we can
    extend $\gamma$ to the automorphism of the splitting field of $f\Phi_j$
    by putting $\gamma(\xi) = \xi^{-1}$. Further let us consider any
    extension of $\gamma$ to an automorphism of $H_{n,j}$,
    which we keep denoting by $\gamma$.

    Let $z$ be a zero
    of the binomial $x^{r^n}-\xi$. It is clear that $\gamma(z)$ is
    a zero of the binomial $x^{r^n}-\xi^{-1}$, i.e.,
    $\gamma(z) = z^{-1}\zeta_n^s$ for some $s$ between $0$ and $r^{n}-1$.
Since by the proof of Theorem 10 $H_{n,j}$ is the radical
extension of the splitting field of $f\Phi_j$, any value of $s$
is possible and it determines completely the automorphism $\gamma$.

As discussed before
    $y_k =z\zeta_n^k+z^{-1}\zeta_n^{-k}$, for $k$ between $0$ and $r^{n}-1$,
      are all the different $r^n$ zeroes of the DeMoivre's polynomial
      $g_n$. We have that $\gamma(y_0) =z^{-1}\zeta_n^{s}
      +z\zeta_n^{-s} =y_{r^n-s}$. Hence the only $\gamma$ which
      fixes all the zeroes of $g_n$ is the one with $s=0$, i.e.,
      $\gamma(z) = z^{-1}$.

      Thus we have established
      that there is exactly one automorphism $\gamma$ of $H_{n,j}$ which
      splits the zeroes of $f\Phi_j$ into $2$-cycles, and fixes
      the zeroes of $g_n$. Hence by the Frobenius-Chebotarev
      Density Theorem the set of primes with such
      a  Frobenius substitution has the prime density equal to
      $\left[H_{n,j}:\mathbb Q\right]^{-1} = \left(2(r-1)r^{n+j-1}\right)^{-1}$.

      We conclude by Theorem 7 that the set of
      primes $\mathcal M_{n}\cap\mathcal K_{j},
      n \leq j$, has the prime density $\left((r-1)r^{n+j-1}\right)^{-1}$.

      To get the prime density of $\Pi_j(t)$ we need to consult
      the Table 1. The prime density 
      \[
      |A_1|=
      |\mathcal K_1| - |\mathcal K_2| -
      |\mathcal M_1\cap \mathcal K_1| + |\mathcal M_1\cap \mathcal K_2|
      = \frac{r-1}{r^2}.
      \]
      Similarly $|A_2|      = \frac{r-1}{r^4}$, and in general $|A_n| =
      \frac{r-1}{r^{2n}}$. It follows by Proposition 4 that the prime
      density $|\Pi_1(t)| =  \frac{1}{r+1}$.

      In the same way for $j\geq 2$ by Proposition 3
      \[
      |C_1\cap \Pi_j(t)|= |\mathcal K_j| - |\mathcal M_1\cap \mathcal K_j| -
      |\mathcal K_{j+1}| + |\mathcal M_1\cap \mathcal K_{j+1}|      =
      \frac{r-1}{r^{j+1}}.
      \]
      Similarly     for $j\geq 2$ and $n\geq 1$
      the prime density $|C_n\cap \Pi_{j}(t)|= \frac{r-1}{r^{2n+j-1}}$.
      It follows by Proposition 4 that for $j\geq 2$ the prime density
      $|\Pi_j(t)| =  \frac{1}{(r+1)r^{j-1}}$.
      The proof of Theorem 2 is completed.
         \end{proof}
  It is instructive to formulate another description of
  the densities of all the sets in the last proof. We introduce
  a shorthand term: a sequence of numbers $\{a_k\}, k\geq k_0$,
  is {\it $r$-geometric} if $a_{k+1} = a_k/r$. We have that for a fixed $n\geq 0$
  the sequences $|\mathcal K_j\cap \mathcal M_n|, j\geq n+1$,
  and $|\mathcal K_j\cap \mathcal M_{n+1}|, j\geq n+1$, are 
  $r$-geometric. It follows that the sequences
  $|\mathcal K_j\cap
  \left( \mathcal M_n\setminus \mathcal M_{n+1}\right)|, j\geq n+1$,
and $|\left( \mathcal K_j\setminus \mathcal K_{j+1}\right)\cap 
\left( \mathcal M_n\setminus \mathcal M_{n+1}\right)|, j\geq n+1$,
are also $r$-geometric. Further, the last sequence of sets coincides with
$\Pi_{j-n}\cap \left( D_n\setminus D_{n+1}\right), j\geq n+1$.
Finally, it leads to
the conclusion that the sequence $|\Pi_k|, k\geq 1$ is $r$-geometric,
i.e., $|\Pi_k| = const/r^k$.

\section{Non-generic values of $t$}

 We address first the non-primitivity of $t$ for a prime $r\geq 2$.
\begin{proposition}
  If a rational $t\neq 0,\pm 1, \pm 2$ then  $t_r = C_r(t)
  \neq 0,\pm 1, \pm 2$ and 
  for any prime $p\in \Pi(t_r)=\Pi(t)$ not dividing the numerator
  of $U_r(t)$, if $r\nmid \chi(t, p)$ then
  $ \chi(t_r, p) = \chi(t, p)$, and if  $r\mid \chi(t, p)$ then
   $\chi(t_r, p) = \chi(t, p)/r$.   Further 
  \[
  \Pi_0(t_r) = \Pi_0(t)\cup \Pi_1(t),
  \ \  \Pi_j(t_r) = \Pi_{j+1}(t), \ \text{for} \ j= 1,2,\dots .
  \]
  \end{proposition}
\begin{proof}
  By the results of \cite{L-W} there is a canonical ring isomorphism
  $\Psi_r :\mathcal R(t) \to \mathcal R(t_r)$ by the conjugation
  $\Psi_r(X) = A^{-1}XA$, for every $X$,
  where  $A = \left[ \begin{array}{cc} 0 &     -U_{r-1}(t) \\
    1  &   U_r(t) \end{array} \right]$.
  Moreover we have $\Psi_r(D_t^r)= D_{t_r}$. To pass to the respective
  isomorphism of $\mathcal S_p(t)$ and  $\mathcal S_p(t_r)$ we need to
  exclude the primes $p$ such that $U_r(t) = 0 \mod p$,
  which by necessity belong to $\Pi_1(t)$.
  By the force of the identity $C_r^2(t) -4 = (t^2-4)U_r^2(t)$ such
  primes belong  also to $\Pi_0(t_r)$
  (remember that we excluded  $r$ from any $\Pi(t)$).
  This information is compatible
  with the last claims of the Proposition.

  For the rest of the primes, identifying $\mathcal S_p(t)$ and
  $\mathcal S_p(t_r)$ by the isomorphism
  we have that $D_{t_r} = D_t^r$, which gives us immediately the formula for
  $\chi(t_r,p)$ in terms of  $\chi(t,p)$. The relation of  the partitions
  follows.
\end{proof}
By this Proposition without loss of generality we can restrict our attention
to $r$-primitive values of $t$.

Let us recall that the $r$-genericity condition appears first
in the proof of Theorem 10. Hence we need to explore how to modify it  for
the non-generic $r$-primitive values of $t$.
\begin{proposition}
   For a rational $t$ if $r =1 \mod 4$ and $t^2-4 = rb^2$ for a rational
  $b\neq 0$ then for any
  prime $p\in \mathcal K_1(t)$ we have
  $r \mod p$ is a square and $ p = 1 \mod r$.
\end{proposition}
In particular in view of Theorem 5 it means that for such $t$
and $p\in \mathcal K_1(t)$ the polynomial $f\Phi_1\mod p$ does not split
quadratically. As a result Theorem 7 acquires a new form
that for any natural $n$ and $j\geq n$ the set of primes $\mathcal M_n\cap
\mathcal K_j$ contains exactly the set of primes $p\in \Pi(t)$ such that
the polynomial $g_nf\Phi_j \mod p$ splits linearly. Moreover since
now in Theorem 10 $\xi \in \mathbb Q(\zeta_j)$ the degree of the extension
$[H_{n,j}:\mathbb Q] = (r-1)r^{n+j-1}$. These two facts lead to the same values
for the prime densities in Theorem 2.
\begin{proof}{(of Proposition 12)}
  For $p\in\mathcal K_1$ by Theorem 1 either $p = 1 \mod r$ or $p=-1 \mod r$.
  Since  $r = 1 \mod 4$ we get that   if $p\in\mathcal K_1$ then
$p\mod r$ is a square.
  We apply now the Gauss law of quadratic reciprocity
  to conclude that for $p\in\mathcal K_1$ also $r \mod p$ is a square.
  By Theorem 1 we finally get that since $t^2-4 =ra^2 \mod p$ is a square
  then $p = 1 \mod r$. 
  \end{proof}
Proposition 12 can also be proven by inspecting the Galois group
of the cyclotomic field of order $r$, i.e., the field $\mathbb Q(\zeta_1)$.
The only automorphism of this field which splits the zeroes of
$\Phi_1(x) = \frac{x^r-1}{x-1}$ into $2$-cycles is the complex conjugation.
Since $\sqrt{t^2-4} = b\sqrt{r}$ is real it is fixed by the automorphism,
  and hence there is no $p$ such that $f\Phi_1\mod p$  splits quadratically.

  The story is completely different for $r = 3 \mod 4$ and
  $t^2-4 = -rb^2$. In this case the complex conjugation provides the unique
  automorphism of all the fields $H_{n,j}$ which fix the zeroes of $g_n$ and
  split the zeroes of $f\Phi_j$ into $2$-cycles. It is so because now
  $\bar{\xi} = \xi^{-1}$ and since $-2 < t < 2$
  all the zeroes of the DeMoivre's polynomial $g_n(x) =C_{r^n}(x)-t$ are real.
  The proof of Theorem 2 can be thus repeated with obvious modifications
  except for $r=3$, which needs to be excluded because in the proof
  of Theorem 10 we used the property that $\mathbb Q(\xi)$ does not contain
  an $r$-th root of unity.   Thus assuming that $r\neq 3$ we get that
  the prime density of the set $\mathcal M_n\cap \mathcal K_j$ is
  now equal to $\frac{2}{(r-1)r^{n+j-1}}$. It leads to doubling of all
    the densities calculated at the end of the proof of Theorem 2.
    As a result we obtain the following modification of Theorem 2.
    \begin{theorem}
\label{Theorem13}
For any prime $r = 3 \mod 4, r\geq 7$, and an $r$-primitive $t$
such that $t^2-4 = -rb^2$ for a rational $b\neq 0$, the disjoint
subsets of primes $\Pi_j(t,r), j= 1,2,\dots,$
have the prime densities and they are equal to
\[
|\Pi_j| = \frac{2}{(r+1)r^{j-1}}.
\]
      \end{theorem}
    We will  prove in the next Section that, under an additional condition,
    Theorem 13 holds also
    in the case $r=3$ and $t^2-4 = -3b^2$, for some
    rational $b\neq 0$, when the field
    $\mathbb Q(\xi)$ is the cyclotomic field of order $3$.
    It was studied by Ballot, \cite{B3}, who calculated the density 
    $|\Pi_0(t)| = 1/4$, under some additional conditions, which is consistent
    with the formula in Theorem 13 for $r=3$. We will shed more
    light on this case, called {\it cubic} in \cite{L-W}.
    First we recall the information from that paper. The crucial fact
    is that in the cubic case there are   roots of unity 
    of order $3$ in the
    special sequence group $\mathcal S(t)$,
    \[
    S = \frac{1}{2b}[ 2 \ 2+b], \ R=   \frac{1}{2b}[ 2 \ 2-b],
    \ S^3=I=R^3, R= S^2=S^{-1}.
    \]
    For any cubic value
    $t$ there are two other rational numbers $a_1 = \frac{-t+3b}{2}$
    and $a_2 = \frac{-t-3b}{2}$, called
    {\it associates}, which satisfy
  $ A=: C_3(t)= C_3(a_1) =     C_3(a_2)$,
    where $C_3(x) = x^3-3x$. We denote $t =a_0$ to stress
    the  symmetry between the three values $a_0,a_1,a_2$. In particular
    all three, and $A$, are cubic, each for its own value of $b$.  
    By the proof of Proposition 11 we get isomorphisms of 
the four rings of matrices $\mathcal R(a_j), j= 0,1,2$, and
$\mathcal R(A)$. Identifying the four
special sequence groups with $\mathcal S(t)$ we have that $D_A = D_t^3$ and 
    $D_{a_1} =SD_t, \ D_{a_2} = RD_t$.
    \begin{proposition}
      For a cubic $t=a_0$  and its associate values $a_1,a_2$
      the sets $     \Pi_0(a_0), \Pi_0(a_1), \Pi_0(a_2) $ are
      mutually disjoint and for any three different indices $k_1, k_2, k_3$
      from $\{0,1,2\}$
      \begin{equation}
        \Pi_1(a_{k_1})= \Pi_0(a_{k_2}) \cup  \Pi_0(a_{k_3}).
        \end{equation} 
      Moreover \begin{equation}
     \Pi_j(a_0)= \Pi_j(a_1)= \Pi_j(a_2), \ \text{for} \ j\geq 2. 
      \end{equation}
    \end{proposition}
    Note that $t$ and $b$ as reduced fractions have the same denominator,
    so that $\Pi(t) = \Pi(a_1)=\Pi(a_2)$. The reason that we have not
    excluded the divisors of $\delta = t^2-4$ from $\Pi(t)$ is exactly that
    we could make this claim. Otherwise the formulation of this Proposition
    would be more complicated. 
    \begin{proof}
      By Proposition 11 we get
      \[
      \Pi_0(A)= \Pi_0(a_{0})\cup \Pi_1(a_{0})=
      \Pi_0(a_{1})\cup \Pi_1(a_{1})=
      \Pi_0(a_{2})\cup \Pi_1(a_{2}),
      \]
      and the equality of the sets in (4).

      We begin with the inspection of primes dividing the
      discriminant $A^2-4$. There are only finitely many of them,
      but it is a union of the primes dividing one of the discriminants
      $a_k^2-4. k=0,1,2$. Indeed, we have
      $A^2-4 = (t^2-4)U_3^2(t) = (a_1^2-4)U_3^2(a_1)$. Hence
      any prime $p$ dividing $A^2-4$ either divides $t^2-4$ and
      belongs to $\Pi_0(t)$, or it divides       $U_3(t) =t^2-1$
      and then the index $\chi(t,p)$ is equal to $3$ or $6$,
      so that $p\in \Pi_1(t)$.

      Since $t^2-4 = -3b^2$ and 
      $a_1^2-4 = \frac{-3(t+b)^2}{4}$, they cannot have
      a common divisor $p\geq 5$. It follows that if $p$ divides $t^2-4$ then
      it belongs to $\Pi_0(t) \cap  \Pi_1(a_1)$.
      Further it is impossible for any prime $p\geq 5$
      that  $t^2= a_1^2 =a_2^2 =  1 \mod p$, so that if $p$ divides $U_3(t)$,
      and hence       $p \in \Pi_{1}(t)$, then       $p \in \Pi_{0}(a_1)$ or
      $p \in \Pi_{0}(a_2)$.  Taking into account the symmetric role played
      by $a_0, a_1,a_2$ we get that divisors of $A^2-4$ satisfy
      all the claims in the Proposition. For the rest of  primes  we can use
      the isomorphism of $\mathcal S_p(a_0)$,
      $\mathcal S_p(a_1)$ and $\mathcal S_p(a_2)$

For $p \in \Pi_0(t)$ and $\chi=\chi(t,p) = \pm 1 \mod 3$ we, have 
$D_t^{\chi} = I \mod p$ and further $\left(DS\right)^{\chi}$ is equal to
$S$ or $R$. It follows that $\left(DS\right)^{3\chi} =I$ and it is
clear that $\chi(a_1,p) = 3\chi$, and hence  $p \in \Pi_1(a_1)$.
In the same way we establish that $p \in \Pi_1(a_2)$. In particular
$p\notin \Pi_0(a_1)\cup \Pi_0(a_2)$. By the symmetry of the problem we
get that the three sets $\Pi_0(a_k), k=0,1,2,$ are mutually disjoint.
It follows further that
$  \Pi_0(a_{k_2}) \cup  \Pi_0(a_{k_3}) \subset \Pi_1(a_{k_1})$
  for any three different indices $k_1,k_2,k_3$. To obtain the opposite
  inclusion let us consider a prime $p \in \Pi_1(t)$
  and let $\chi' = \chi/3 = \pm 1 \mod 3$. We have that
  $D^{\chi'} = S\ \text{or} \ R$ and accordingly $(DS)^{\chi'}$ or
  $(DR)^{\chi'}$ is equal to the identity $I$, which means that 
  $ \Pi_1(t)\subset  \Pi_0(a_{1}) \cup  \Pi_0(a_{2}) $.
  The symmetry delivers the rest.
      \end{proof}

    \section{Cubic case and cubic primitivity}

    \begin{definition}
      A rational $t$ is called cubic if $t^2-4 = -3b^2$ for
      a rational $b\neq 0$. It  is called {\it cubic primitive}
      if $t$ and its associate values $a_1$ and $a_2$ are all $3$-primitive.
      \end{definition}

    \begin{proposition}
      If $t$ is cubic and cubic primitive then
      $\mathbb Q(\zeta_1,z_1)\cap\mathbb Q(\zeta_2)=
      \mathbb Q(\zeta_1)$ and
      $Gal(\mathbb Q(\zeta_2,z_1)/\mathbb Q(\zeta_1))$ is isomorphic
      to $\mathbb Z_3 \oplus \mathbb Z_3$.
    \end{proposition}
    \begin{proof}    
Since 
$[\mathbb Q(\zeta_1,z_1):
      \mathbb Q(\zeta_1)]=3$ and
      $[\mathbb Q(\zeta_2):
      \mathbb Q(\zeta_1)]=3$, to get the first claim
      we just need to exclude the possibility that
      $\mathbb Q(\zeta_1,z_1)=\mathbb Q(\zeta_2)$.
      These fields are, respectively, the splitting fields of
      $z^3-\xi$ and $z^3-\zeta_1$ over $\mathbb Q(\zeta_1)$.
      It follows from the theorem of Schinzel
      (\cite{S1}, also \cite{B3}, page 82),
      that such fields coincide if and only if
      $\xi = \zeta_1^{\epsilon} \beta^3$ for $\beta\in \mathbb Q(\zeta_1)$
      and $\epsilon =\pm 1$.

      For a cubic $t$ we have that
      $\mathbb Q(\xi)=\mathbb Q(\zeta_1)$ and we can identify this field
      with $\mathcal R(t)$, so that $\xi$ and $\zeta_1^{\epsilon}$ become
      $D$ and $S$ or $R$, respectively.  The condition for the equality of the
      splitting fields becomes then $D=SB^3$ or $D=RB^3$
      for some $B\in \mathcal R(t)$.
      It follows that $1 =\det D = (\det B)^3$ and by necessity
      $B\in \mathcal S$. By definition $a=tr \ B^3$ is an associate value for
      $t$ and it is not $3$-primitive. The contradiction completes
      the proof of the
      first part of the Proposition. The second part follows from the
      elementary Galois Theory.
\end{proof}
    In particular for a cubic primitive $t$ the binomial
    $z^3-\xi$ is irreducible over $\mathbb Q(\zeta_2)$
    and the Galois group $Gal(\mathbb Q(\zeta_2,z_1)/\mathbb Q(\zeta_1))$
    is not cyclic.
    \begin{corollary}
      For a cubic primitive $t$ the binomial $z^{3^n}-\xi$ is
      irreducible over $\mathbb Q(\zeta_j)$, for $j\geq n = 1,2,\dots $. 
    \end{corollary}
    \begin{proof}
      We only need to consider $n\geq 2$.
      Let, as before, $z_n$ be a root of the binomial  $z^{3^n}-\xi$.
      We have that  $z_1=z_n^{3^{n-1}}$ is a zero of the binomial $z^3-\xi$.
      If $z_n$ belongs to the cyclotomic field $\mathbb Q(\zeta_j)$ then
        the splitting field of $z^3-\xi$ over $\mathbb Q(\zeta_2)$
	would be a subfield of $\mathbb Q(\zeta_j)$. It is a contradiction
	since the larger field has a cyclic Galois group over
	$\mathbb Q(\zeta_1)$
	and the smaller field
        a non-cyclic Galois group over 	$\mathbb Q(\zeta_1)$.
    \end{proof}   
    With the help of Corollary 16 we can claim that in Theorem 10
    $[\mathbb Q(\zeta_j,z_n):\mathbb Q(\zeta_j)] = 3^n$, for $j\geq n\geq 1$
    also for a cubic $t$ in the cubic primitive case.    
    As a result Theorem 13 is valid also for $r= 3$,
    if the cubic primitivity is added to the assumptions.

    Finally let us address the case when  cubic primitivity is not
    satisfied. Towards that end let 
    $t' = C_{3^k}(t)$ for a fixed $k\geq 1$ and a
    cubic primitive $t$, and let $a_1$ and $a_2$ be
    the associate values of $t'$, which are clearly not cubic primitive. 
    We are going to describe the partition for $a_1$ in terms of
    the partition for $t$.
    Proposition 11 gives us the connection between the partitions for $t'$ and
    $t$. 
    Using  Propositions  14 we obtain further 
    \begin{equation}
    \begin{aligned}
      &\Pi_0(a_1)\cup \Pi_0(a_2) =\Pi_1(t')= \Pi_{k+1}(t),  \\
      &\Pi_1(a_{1}) = \Pi_0(a_2)\cup \Pi_0(t') = \Pi_0(a_2) \cup  
      \bigcup_{j= 0}^k\Pi_{j}(t).
      \end{aligned}
  \end{equation}
   and $\Pi_j(a_1)=\Pi_j(t') = \Pi_{j+k}(t)$ for $j=2,3,\dots$.
    To obtain the densities for the partition for $a_1$ in terms of
    the densities for the partition for $t$ it remains to establish
    how $\Pi_{k+1}(t)$ is split into two
    disjoint subsets $\Pi_0(a_1)$ and $\Pi_0(a_2)$.

    We claim that the densities of $\Pi_0(a_1)$ and $\Pi_0(a_2)$
    must be equal. Identifying the fields $\mathcal R(a_1)$
    and $\mathcal R(a_2)$ with $\mathcal R(t)$ we get
    $D_{a_1}  = S(D_{t})^{3^k}$ and     $D_{a_2}  = R(D_{t})^{3^k}$.
    It follows that $a_1$ and $a_2$ are by necessity $3$-primitive.
    With the further identification
    of $\mathcal R(t)$ with the cyclotomic field $\mathbb Q(\zeta_1)$ 
    we obtain that $\xi_{a_1} = \zeta_1\left(\xi_{t}\right)^{3^k}$ and
    $\xi_{a_2} = \zeta_1^{-1}\left(\xi_{t}\right)^{3^k}$.
    Let as before $z_n$ be a root of the binomial $z^{3^n}-\xi_t$.  
    We will argue that the fields
    $H_{n,j}(a_1)=\mathbb Q\left(\zeta_j,\xi_{a_1},
    \zeta_{n+1}z_n^{3^k}\right)=\mathbb Q\left(\zeta_j,\xi_{t},
    \zeta_{n+1}z_n^{3^k}\right)$ and
    $H_{n,j}(a_2)=\mathbb Q\left(\zeta_j,\xi_{a_2},
    \zeta_{n+1}^{-1}z_n^{3^k}\right)=\mathbb Q\left(\zeta_j,\xi_{t},
    \zeta_{n+1}^{-1}z_n^{3^k}\right)$
    are extensions of the same degree over $\mathbb Q$
    for $j\geq n\geq 1$.
    Actually it is clear that
    $H_{n,j}(a_1)
    =\mathbb Q\left(\zeta_j,\zeta_{n+1},\xi_{t}\right) =
    H_{n,j}(a_2)$
    for $n \leq k$. Also obviously  $H_{n,j}(a_1) = H_{n-k,j}(t)=
    H_{n,j}(a_2)$ for $j\geq n+1, n\geq k+1$.
    We need to inspect further only the case $j=n\geq k+1$.
    For such $n$  
    we have $[H_{n,n+1}(a_1):H_{n-1,n}(a_1)] =[H_{n-k,n+1}(t):H_{n-k-1,n}(t)] =9$
    because of the assumption that
    $t$ is cubic primitive. In the tower of extensions
    $H_{n-1,n} < H_{n,n} < H_{n,n+1}$
    the intermediate extensions have at most degree $3$, and hence they must
    have both  degree $3$.

    We conclude that the  densities of the  partitions for $a_1$ and
    $a_2$ must coincide.

    \section{$r=2$}

    Theorem 5 is not true as stated for $r =2$. It can be modified
    to include this case by invoking explicitly the
    permutation  $\sigma$ of zeroes of $f\Phi_j$, which maps any
    zero $\zeta$ of $\Phi_j$ into its inverse $\zeta^{-1}$,
    and also $\xi$ into $\xi^{-1}$.
    The resulting Theorem 5 states that $p \in \mathcal K_j$ iff
    $f\Phi_j \mod p$ splits linearly or the Frobenius substitution 
    $\sigma_p = \sigma$. With enough genericity conditions it is then
    possible to establish the validity of Theorem 2 also for $r=2$.
    The trouble with this approach is the difficulty in complete enumeration of
    non-generic cases. It was achieved by a different scheme in \cite{W},
    and it is quite involved. We will obtain it now in a different
    way by modifying Theorem 5 into the following claim.
    Let $\widetilde f(x) = x^2 +t^2-4$ and $c_j(x) = C_{2^j}(x)$,
    where $C_n(x)$ is the Chebyshev polynomial of the first kind.
    Let us note that for $r=2$ the set $\mathcal K_1(t) = \Pi(t)$,
    because the special sequence group $\mathcal S_p(t)$ always contains
    a square root of the identity $I$ equal to $-I$.
    \begin{theorem}
  For $r =2$ and a natural $j\geq 2$
  the set $\mathcal K_j(t)$
  contains exactly the primes
  $p\in \Pi(t)$  such that
  the polynomial $\widetilde fc_{j-2}\mod p$ splits linearly.
    \end{theorem}
    \begin{proof}
      Let us start with $j=2$. We want to characterize the primes $p$
      such that the special sequence group $\mathcal S_p(t)$ contains
      a square root $J$ of $-I$. Since $J^2 =-I$ we conclude that for
      $a = tr\ J$ we have 
      $-2 = tr\ J^2 = C_2(a) = a^2 -2$, and thus $a=0$.
      Using Lemma 6 we obtain that $\mathcal S_p(t)$ contains the
      matrix $J$ with trace $0$ if and only if $-(t^2-4)$ is a square
      $\mod p$, or equivalently $\widetilde f(x) \mod p$ splits.
      The condition on $c_0(x) = C_1(x) =x$ is void.

      In general for $j\geq 3$ we want to characterize the primes $p$
      such that the special sequence group $\mathcal S_p(t)$ contains
      a matrix $J$ with zero trace and also a  root of
      $J$ of degree $2^{j-2}$. The last condition is by Lemma 8 equivalent to
      $c_{j-2}(x) \mod p$ having a zero $b$. The respective matrix
      $B \in \mathcal S_p(t)$ with the trace equal to $b$
      is a primitive root of unity of order
      $2^j$. The traces of powers of $B$ give us all the other roots
      of $c_{j-2} \mod p$, hence it splits linearly.

      The remaining logical considerations are routine and we leave them
      to the reader.
    \end{proof}
    The proof of Theorem 7 is also valid for $r=2$ except for the
    use of Theorem 5. We can now reformulate Theorem 7 using Theorem 17.
    \begin{theorem}
  For  $r=2$ and any natural $n$ and $j\geq n, j\geq 2$
  the set $\mathcal M_n\cap \mathcal K_{j}$
  contains exactly the primes
  $p\in \Pi(t)$  such that the polynomial $c_{j-2}\widetilde fg_n\mod p$
  splits linearly.
      \end{theorem}
    
\begin{definition}
    A rational $t\neq 0,\pm 1, \pm 2$ is {\it twin primitive}
    if both $t$ and $-t$ are $2$-primitive.

    A rational $t\neq 0,\pm 1, \pm 2$ is $2$-generic
    if 
    $[\mathbb Q(\sqrt{2+t},\sqrt{2-t},\sqrt{2}):\mathbb Q]=8$.
\end{definition}
Twin primitivity  is equivalent to both
$2+t$ and $2-t$ not being rational squares.
The role of twin primitivity 
will be illuminated in the next section. Clearly, if a rational
$t$ is $2$-generic then it is twin primitive.

We denote by  $\widetilde H_{n,j}$ the splitting field of the
polynomial $c_{j-2}\widetilde fg_n$. Further for $j\geq 3$  let $L_j$ 
be the splitting field of the polynomial $c_{j-2}= C_{2^{j-2}}$.

Since $C_{2^k} = (C_2)^k = C_2\circ C_2 \circ ... \circ C_2 $,
the zeroes of $c_{k}$ and $g_k=c_{k}-t$ can be
expressed by the nested quadratic roots
\[
b_k = \pm\sqrt{2\pm\sqrt{2\pm\dots \sqrt{2}}},\ \ \ \text{and} \ \ \
y_k=\pm\sqrt{2\pm\sqrt{2\pm\dots \sqrt{2+t}}}.
\]
\begin{proposition}
  For $j\geq 3$ the field $L_j$ is the real part of the
  cyclotomic field of order $2^j$. Moreover
  $L_j = \mathbb Q\left(b_{j}\right)$, where $b_j$ is
  any root of $c_{j-2}$ and
  $Gal\left(L_j/\mathbb Q\right)$ is cyclic of order $2^{j-2}$.
\end{proposition}
This proposition is part of math folklore. We provide a proof for
the convenience of the reader.
  \begin{proof}
    We will rely on the following formula for powers of complex numbers.
    For  $z= x +iy, |z| = 1$ we have
    \[
    z^k = \frac{1}{2}\left( C_k(2x) + i y U_k(2x)\right).
    \]
    It follows  that  $z_0=\frac{1}{2}\left(b_j + i\sqrt{4-b_j^2}\right)$
    is a primitive root of unity of order $2^j$.
    Indeed, since $C_{2^{j-2}}(b_j) =0$ then  using  (2)
    we get $\sqrt{4-b_j^2} = 2/U_{2^{j-2}}(b_j)$,
    so that  $z_0^{2^{j-2}} = i$. 
    Further,
      we can obtain any root of unity of  order $2^j$ by taking powers
      of $z_0$. Their real parts belong to $\mathbb Q\left(b_{j}\right)$,
      which is thus equal to $L_j$ and also
      the real part of the cyclotomic field.

      The claim about the order of the Galois group is obtained  by
      elementary Galois theory.
\end{proof}

  We proceed to  find the degrees  of the field extensions
  $\widetilde H_{0,j}, j\geq 2$ over $\mathbb Q$. 

  \begin{proposition}
        For any rational $t \neq 0,\pm 1, \pm 2$ we have  $\widetilde H_{0,2}=
    \mathbb Q\left(\sqrt{4-t^2}\right)$ and for $j\geq 3$ .
    \[
    \widetilde H_{0,j}=\mathbb Q\left(\sqrt{4-t^2},b_{j}\right) =
    L_j\left(\sqrt{4-t^2}\right).
    \]
    For any $2$-generic $t$ and $j\geq 2$
     \[
    [\widetilde H_{0,j}:\mathbb Q] = 2^{j-1}.
    \]
    \end{proposition}
  \begin{proof}
    The first part follows directly from definitions.
    To prove the second part we use the fact that $L_j$, having
    a cyclic Galois group, contains only one quadratic subfield 
    $\mathbb Q\left(\sqrt{2}\right)$. By the genericity assumption
    about $t$
    it follows  that $\sqrt{4-t^2}$ does not belong to $L_j, j\geq 3$,
    and hence
    $[L_j\left(\sqrt{4-t^2}\right):L_j] =2$.
    \end{proof}
  We proceed with further exploration of the fields
  $\widetilde H_{n,j}$. We have
  $\widetilde H_{1,2}=\mathbb Q\left(\sqrt{4-t^2},\sqrt{2+t}\right)=
  \mathbb Q\left(\sqrt{2+t},\sqrt{2-t}\right)$ and
  $\widetilde H_{1,3}=
  \mathbb Q\left(\sqrt{2+t},\sqrt{2-t},\sqrt{2}\right)$.
  For any $2$-generic $t$ these field extensions over $\mathbb Q$
  have orders $4$ and $8$ respectively.
  \begin{theorem}
    For any rational $t \neq 0,\pm 1, \pm 2$ and any $n\geq 1$ and
    $j\geq n, j\geq 3,$
    \[
\widetilde H_{n,j}=L_j\left(\sqrt{4-t^2},y_{n}\right),
\]
for any root $y_n$ of the polynomial $g_n(x) = c_n(x)-t$.

Further $[\widetilde H_{n,j}:\mathbb Q] = 2^{n+j-1}$, for any $2$-generic $t$. 
  \end{theorem}  
  \begin{proof}
    Our first task is to establish that
    $L_j\left(\sqrt{4-t^2},y_{n}\right)$ contains all the roots of $g_n$.
    To that end we consider the ring $ \mathcal R(t,K)$ with
    $K=L_j\left(\sqrt{4-t^2},y_{n}\right)$. By Lemma 8 the ring contains
    a primitive root of unity of degree $2^j$ and
    the matrix $Y$ with determinant $1$ and trace $y_n$, which
    is a root of degree $2^n$ of $D$. Since $j\geq n$ the ring contains also
    $2^n$ matrices which are different roots of unity of degree $2^n$.
    Multiplying $Y$ by these roots of unity
    we obtain matrices with traces equal to all the other zeroes of $g_n$.
    This sketchy argument can be completed in the same way as in the
    proof of Theorem 7, which we leave to the reader.

    Since $y_{n} = \sqrt{2+y_{n-1}}$ we conclude that for 
    $j\geq n\geq 2,$
    $[\widetilde H_{n,j}: \widetilde H_{n-1,j}]$ 
    is either $1$ or $2$.
    Our goal is to establish that for any $2$-generic $t$
    the order of extension
    is always $2$. To that end we will employ an elementary technical
    lemma.

    Let $K_0$ be a field of characteristic $0$ and $a_0\in \mathbb Q,
     a_0\neq \pm 2$.
  We define $K_1= K_0(\sqrt{4-a_0^2})$,
  and $K_2= K_1(a_1)$ where $a_1= \sqrt{2+a_0}$. We have that
  $\widetilde a_1 =\sqrt{2-a_0} \in K_2$.
  Further, by induction, given a field
  $K_s=K_{s-1}(a_{s-1}), a_{s-1}=\sqrt{2+a_{s-2}},
  \widetilde a_{s-1}=\sqrt{2-a_{s-2}},
  s\geq 2$, we define
  $K_{s+1}=K_s(a_s), a_s=\sqrt{2+a_{s-1}}$ and   $\widetilde a_s
  =\sqrt{2-a_{s-1}}$.

  Let us observe first that $a_{s}\widetilde a_{s} =
  \widetilde a_{s-1} \neq 0$. Indeed we have  
  \[
    a_{s}\widetilde a_{s} =\sqrt{2+a_{s-1}} \sqrt{2-a_{s-1}}=\sqrt{4-a_{s-1}^2}
    = \sqrt{2-a_{s-2}} =  \widetilde a_{s-1}.
    \]
  It follows by induction that     $a_{s}\widetilde a_{s}\neq 0$ and 
  $\widetilde a_{s}\in K_{s+1}=K_s(a_s)$. Consequently
  $K_s(\widetilde a_s)=K_s(a_s)$.
\begin{lemma}
    If $[K_{s+1}:K_s] =2$, for some $s\geq 1$,
    then     $[K_{r+1}:K_{r}] =2$ for every $r\geq s$.
  \end{lemma}
  \begin{proof}
    Suppose to the contrary that $m\geq 1$ is the first natural number such
    that $K_{m+2}=K_{m+1}$. It means that $a_{m+1}=\sqrt{2+a_m}
    \in K_{m+1}= K_m(a_m)$.
    Hence there are $\alpha,\beta \in K_m$ such that
    $2+a_m = (\alpha a_m +\beta)^2$. It follows that 
    $\alpha^2a^2_m+\beta^2 =2$, $2\alpha\beta =1$.
    Note that while $a_m\notin K_m$ by the inductive
    assumption, $a_m^2\in K_m$     by definition.
    Combining the two equations we arrive at
    $4(\beta^2-1)^2 = 4-a_m^2$. Hence $\sqrt{2-a_{m-1}} 
    =\sqrt{4-a_m^2} = 2(\beta^2-1)$. We arrived at a contradiction
    with the inductive assumption that $\sqrt{2+a_{m-1}}$, and hence also 
    $\sqrt{2-a_{m-1}}$, are not in the field $K_m$. It proves that
    actually $\sqrt{2+a_m} \notin K_{m+1}$.
  \end{proof}

  To finish the proof of the theorem we apply the lemma with
  $K_0= L_j$ and $a_0 =t$. We get that $K_1 = L_j\left(\sqrt{4-t^2}\right)$
  and $K_2 = L_j\left(\sqrt{2+t},\sqrt{2-t}\right)$. The subfields of
  quadratic irrationalities in $K_1$ and $K_2$ are
  $\mathbb Q\left(\sqrt{2},\sqrt{4-t^2}\right)$ and
  $\mathbb Q\left(\sqrt{2},\sqrt{2+t},\sqrt{2-t}\right)$,
  respectively. We conclude, using the $2$-genericity of $t$, that
  $[K_2:K_1] =2$.    
\end{proof}

  Using the orders of extensions from Theorem 21 we can now
  repeat the calculations
  of prime densities at the end of proof of Theorem 2 from Section 4.
  In this way we obtain Theorem 2 also in the generic case of $r=2$.

    \section{$r=2$, non-generic cases}

    The case of $r=2$ comes with the {\it twin symmetry},
    $t\leftrightarrow -t$.
    
    \begin{proposition}
      For any rational $t\neq 0,\pm 1, \pm 2$, and a prime $p \in
      \Pi(t) =\Pi(-t)$
      \[
      \Pi_0(-t) = \Pi_1(t), \ \ \ \ \Pi_1(-t) = \Pi_0(t).
      \]
      Further $\Pi_j(-t) = \Pi_j(t)$, for any $j\geq 2$. 
      \end{proposition}
\begin{proof}
  Let us observe that $D_{-t} = -D_t^T$. It follows that for $p\in \Pi_0(t)$
  \[
  D_{-t}^{\chi(t,p)} =   (-1)^{\chi(t,p)}   \left(D_{t}^T\right)^{\chi(t,p)} = -I.
  \]
  Hence $\chi(-t,p) = 2\chi(t,p)$, for $p \in \Pi_0(t)$. We conclude that
  $\Pi_0(t) \subset \Pi_1(-t)$, and by symmetry also
  .$\Pi_0(-t) \subset \Pi_1(t)$.

  If  $p\in \Pi_1(t)$ then for the odd  $k= \chi(t,p)/2$ we have 
\[
  D_{-t}^{k} =   (-1)^{k}   \left(D_{t}^T\right)^{k} = I.
  \]
  Hence $\chi(-t,p) = \chi(t,p)/2$, for $p \in \Pi_1(t)$. We conclude that
  $\Pi_1(t) \subset \Pi_0(-t)$. Taking into account the set inclusion
  above we get $\Pi_1(t) =\Pi_0(-t)$, and by symmetry
  also $\Pi_1(-t) =\Pi_0(t)$.

  Further if $\chi(t,p) = 0 \mod 4$ then for the even $k= \chi(t,p)/2$ 
\[
  D_{-t}^{k} =   (-1)^{k}   \left(D_{t}^T\right)^{k} = -I.
  \]
  It follows that $\chi(-t,p) = 2k= \chi(t,p)$. Hence $\Pi_j(-t) = \Pi_j(t)$
  for any $j\geq 2$.
\end{proof}

Combining the last Proposition with Proposition 11 we get the following
Corollary.
\begin{corollary}
  For any rational $t\neq 0,\pm 1, \pm 2$ and natural $k\geq 1$, if
  $t_{k} = C_{2^k}(t)$
  then 
  \[
  \Pi_0(-t_{k}) = \Pi_{k+1}(t), \ \  \ \
  \Pi_1(-t_{k}) =   \bigcup_{j= 0}^k\Pi_{j}(t),
  \]
  Further,
  $\Pi_j(-t_{j}) = \Pi_{j+k}(t)$,  for any $j\geq 2$.
\end{corollary}
It is striking that the set of primes $\Pi_0(t)$ can  be arbitrarily small
and $\Pi_1(t)$ arbitrarily large, for appropriate values of $t$.


The Corollary makes it clear that it is enough to study the  non-generic
values of $t$ which are twin primitive.
The first non-generic case  that we distinguish is when $4-t^2$ is a
rational square. We will call such a value of $t$   {\it circular}. It can be
characterized by $\mathcal R(t)$ being isomorphic to $\mathbb Q(i)$.
It will be discussed in detail in the next section.

To enumerate the remaining  non-generic cases
we can assume that $\sqrt{2+t}$,
$\sqrt{2-t}$ and $\sqrt{4-t^2}$ are not
rational. We are left with two cases:
non-genericity of {\it type  $A$} when $2(2+t)$, or $2(2-t)$,
is a rational square,
and 
non-genericity of {\it type  $B$} when $2(4-t^2)$ is a rational square.

We start the discussion with  the simpler case of type $B$.
A sequence of numbers $\{a_n\}_{n\in \mathbb N}$
is called {\it $2$-geometric} if 
$a_{n+1} = a_n/2$.
\begin{theorem}
  If  $t$ is non-generic of type $B$ then the densities are 
  \[
  |\Pi_0(t)| =  |\Pi_1(t)| = \frac{7}{24}, \ \ 
    |\Pi_2(t)| = \frac{1}{12}, 
  \]
  and the sequence $|\Pi_k(t)|, k\geq 3$, is $2$-geometric.
  \end{theorem}
\begin{proof}
  If $2(4-t^2)$ is a rational square then
  \[
   \mathbb Q(\sqrt{4-t^2}) =   \mathbb Q(\sqrt{2}) = \widetilde H_{0,2}
   =  \widetilde H_{0,3}.
   \]
   At the same time $ [\widetilde H_{1,3}:  \widetilde H_{0,3}] =
   [\mathbb Q(\sqrt{2},\sqrt{2+t}):  \mathbb Q(\sqrt{2})] = 2$.
   Similarly $[\widetilde H_{1,j}:  \widetilde H_{0,j}] =
   [L_j(\sqrt{2+t}):  L_j] = 2$, for $j\geq 3$. To prove it we can
   repeat the argument from the proof of Theorem 21 that
   $L_j(\sqrt{2+t})\neq  L_j$ because their quadratic subfields
   are $\mathbb Q(\sqrt{2},\sqrt{2+t})$
   and $\mathbb Q(\sqrt{2})$, respectively. It follows now by Lemma 22
   that 
   $[\widetilde H_{n+1,j}:  \widetilde H_{n,j}] =2$ for any
   $n \leq j-1, j\geq 3$.
   
   With this information at hand we can calculate the densities
   of the partitions inside $D_k\setminus D_{k+1}$ for $k\geq 0$, cf. Table 1.
   The general pattern is that $D_0\setminus D_1$ and
   $D_1\setminus D_2$ are special, while for the rest we have the
   same splitting, up to scale,
   as in the generic case, i.e. for a fixed $k\geq 2$
   the sequence
   $|\Pi_j\cap \left( D_k\setminus D_{k+1}\right)|, j\geq 2,$ is $2$-geometric,
   and 
   $|\Pi_0\cap \left( D_k\setminus D_{k+1}\right)|=
   |\Pi_1\cap \left( D_k\setminus D_{k+1}\right)|$.

   It remains to study the partitions inside 
$D_0\setminus D_1$ and
   $D_1\setminus D_2$. Since $\widetilde H_{0,2} = \widetilde H_{0,3}$
   it follows that
   $\mathcal K_3 = \mathcal K_2$. Consequently inspecting Table 1
   we conclude that
   $\Pi_2\cap \left( D_0\setminus D_{1}\right)$ is
   empty and that the sequence
   $|\Pi_j \cap\left( D_0\setminus D_{1}\right)|,j\geq 3,$
   is $2$-geometric. It is also obvious that
   $|\Pi_0\cap \left( D_0\setminus D_{1}\right)|=
   |\Pi_1\cap \left( D_0\setminus D_{1}\right)| = \frac{1}{4}$.

   Further $\Pi_0\cap \left( D_1\setminus D_{2}\right)$
   and $\Pi_1\cap \left( D_1\setminus D_{2}\right)$ are empty
   while the sequence $|\Pi_j \cap\left( D_1\setminus D_{2}\right)|,j\geq 2,$
   is $2$-geometric.

   Taking into account that in this case
   obviously $|D_1| =1/4$ and $|D_2| =1/8$
   we obtain 
   \[
   |\Pi_0| =   |\Pi_1| = \frac{1}{4} + 0 + \frac{1}{8}\cdot\frac{1}{3}=
   \frac{7}{24} \ \ \
   \text{and} \ \ \
   |\Pi_2| =  0+  \frac{1}{16} +  \frac{1}{8}\cdot\frac{1}{6} = \frac{1}{12}.
   \]
  \end{proof}

We proceed with the case of type $A$ when $2(2+t)$ is a rational square.
The case of $2(2-t)$ being a rational square
does not require special attention since the respective partition
can be obtained by Proposition 23 from the partition for $-t$.

We will need the following lemma
about the splitting field of $g_2= c_2(x)-t$.
\begin{lemma}
  If $[\mathbb Q\left(\sqrt{2+t},\sqrt{2-t}\right):
    \mathbb Q] =4$ then the
  Galois group of the splitting field of
  $c_2(x)-t = x^4-4x^2+2 -t = \left(x^2-2\right)^2 -2 -t$
  is non-abelian.
  \end{lemma}
\begin{proof}
  The four roots of  our polynomial are
  $\pm\sqrt{2\pm\sqrt{2+t}}$. Hence the polynomial is irreducible
  and the Galois group acts transitively on these roots.
  Further since $\sqrt{2+\sqrt{2+t}}\sqrt{2-\sqrt{2+t}} =
  \sqrt{2-t}$ we get that $\mathbb Q(\sqrt{2+t},\sqrt{2-t})$
  is a normal subfield of the splitting field with the Galois group
  isomorphic to $\mathbb Z_2\oplus \mathbb Z_2$. Also the splitting field
  is equal to $\mathbb Q(\sqrt{2+t},\sqrt{2-t}, \sqrt{2+\sqrt{2+t}})$
  which is a field extension of degree $8$ over $\mathbb Q$.

  Let us consider the splitting field automorphism $\sigma_1$
  such that $\sigma_1(\sqrt{2\pm t}) = -\sqrt{2\pm t}$.
  By necessity $\sigma_1(\sqrt{2+\sqrt{2+t}}) =
  \pm\sqrt{2-\sqrt{2+t}}$. Let us assume without loss of generality
  that the sign in the last formula is
  plus. We have
  \[
  \sigma_1\left(\sqrt{2+\sqrt{2+t}}\sqrt{2-\sqrt{2+t}}\right)
    = \sigma_1(\sqrt{2-t}) = -\sqrt{2-t}
    \]
    It follows immediately that $\sigma_1$ has order $4$.

    Further we consider the field automorphism $\sigma_2$ such that
    $\sigma_2(\sqrt{2\pm t}) = \mp\sqrt{2\pm t}$.
We have $\sigma_2\left(\sqrt{2+\sqrt{2+t}}\sqrt{2-\sqrt{2+t}}\right)
    = \sigma_2(\sqrt{2-t}) = \sqrt{2-t}$.
    It follows immediately that $\sigma_2$ has order $2$.

    It is somewhat tedious, but elementary, that $\sigma_2\circ\sigma_1
    \neq \sigma_1\circ\sigma_2$.
\end{proof}
For clarity let us note that it can be established that
if we assign $\pm\sqrt{2+\sqrt{2+t}}$ and $\pm\sqrt{2-\sqrt{2+t}}$
to the diagonally opposite vertices of a square then the Galois group
becomes naturally isomorphic to the group of isometries of the square,
i.e., the dihedral group $D_4$.

We are ready to address the non-genericity of type $A$.
\begin{theorem}
  If  $t$ is non-generic of type $A$ then the densities are
  \[
  |\Pi_0(t)|  =  |\Pi_1(t)| = \frac{7}{24}, \ \ 
    |\Pi_2(t)| = \frac{1}{3}, 
  \]
  and the sequence $|\Pi_k(t)|, k\geq 3$, is $2$-geometric.
  \end{theorem}
\begin{proof}
If $2(2+t)$ is a rational square then
  \[
   \mathbb Q(\sqrt{4-t^2},\sqrt{2}) = \widetilde H_{0,3}
   =  \widetilde H_{1,2}=\widetilde H_{1,3}.
   \]
   It implies that $\mathcal K_3 =  \mathcal M_1\cap \mathcal K_2  =
   \mathcal M_1\cap \mathcal K_3$. Hence the densities in
   $D_0\setminus D_1$ (cf. Table 1) are
   \[
   |\Pi_0\cap \left(D_0\setminus D_1\right)| =
   |\Pi_1\cap \left(D_0\setminus D_1\right)| =
   |\Pi_2\cap \left(D_0\setminus D_1\right)| = \frac{1}{4},
\] 
while $\Pi_j\cap \left(D_0\setminus D_1\right) =\emptyset $
for $j\geq 3$.

Further, inside $D_1$  we have
$\Pi_0\cap \left(D_1\setminus D_2\right) =
\Pi_1\cap \left(D_1\setminus D_2\right) = \emptyset $
and the sequence $|\Pi_k\cap  D_1  |, k\geq 2$
is $2$-geometric, because for $j\geq 3$ 
$[\widetilde H_{1,j+1}:\widetilde H_{1,j}] =
 [L_{j+1}(\sqrt{4-t^2}):L_j(\sqrt{4-t^2})] =2$.

 Inside $D_k, k\geq 2$ the situation is, up to scale,
 the same as in the generic case,
 i.e., i.e. for a fixed $k\geq 2$
   the sequence
   $|\Pi_j\cap \left( D_k\setminus D_{k+1}\right)|, j\geq 1,$ is $2$-geometric,
   and 
   $|\Pi_0\cap \left( D_k\setminus D_{k+1}\right)|=
   |\Pi_1\cap \left( D_k\setminus D_{k+1}\right)|$.
 It follows from the fact that 
 $[\widetilde H_{k+1,j}:\widetilde H_{k,j}] =2$ for $k\leq j-1, j\geq 3$.
 To prove it we use  Lemma 22, establishing first that
 for any $j\geq 3$  $[\widetilde H_{2,j}:\widetilde H_{1,j}] =2$.
To prove the last claim we invoke Lemma 26. 
Indeed $\widetilde H_{2,j} \neq \widetilde H_{1,j}$ because
the first field contains the splitting field of $g_2 = c_2(x)-t$
with a non-abelian Galois group, while the second field has
an abelian Galois group.

Using the information obtained so far, and 
taking into account that $|D_1| =1/4$ and $|D_2| =1/8$ we obtain
the following densities.
   \[
   |\Pi_0| =   |\Pi_1| = \frac{1}{4} + 0 + \frac{1}{8}\cdot\frac{1}{3}=
   \frac{7}{24} \ \ \
   \text{and} \ \ \
   |\Pi_2| =  \frac{1}{4}+  \frac{1}{16} +  \frac{1}{8}\cdot\frac{1}{6} =
   \frac{1}{3}.
   \]
\end{proof}

\section{$r=2$, circular case}

If $t$ is circular then there is a rational $w$
such that $t^2+w^2 =4$. This property is parallel
to cubic or twin symmetry.
Clearly such values are in $1-1$ correspondence
with Pythagorean triples. We call $w$ the {\it associate value} for $t$.
\begin{definition}
  A circular $t$ is called {\it circular primitive} if both
  $t$ and its associate value $w$ are $2$-primitive.
  \end{definition}
Let us note that for a circular $t$ the notions of $2$-primitivity
and twin-primitivity coincide. Indeed it follows directly from
$(2+t)(2-t) = w^2$. Further, for a circular $t$ and its associate
value $w$ either both are $2$-primitive or exactly one is, i.e.,
it is impossible that simultaneously $2+t$ and $2+w$ are rational
squares. Indeed we have 
$2(2+t)(2+w) = (2+t+w)^2$, (\cite{L-W}). It also follows from this identity
that $t$ is circular primitive if and only if $2+t$ and $2(2+t)$
are not rational squares, or equivalently
that $[\mathbb Q(\sqrt{2},\sqrt{2+t}):\mathbb Q] =4$.
Let us note that if $t^2+w^2 = 4$ then for any odd $n$ (\cite{W},
\cite{L-W})
\begin{equation}
  C_n(t)^2 +C_n(w)^2 =4.
  \end{equation}
\begin{theorem}
  If  $t$ is circular primitive then 
  \[
  |\Pi_0(t)| =  |\Pi_1(t)| = \frac{1}{6},
  \]
  and the sequence $|\Pi_k(t)|, k\geq 2$, is $2$-geometric.
  \end{theorem}
\begin{proof}
  If $4-t^2$ is a rational square then $\widetilde H_{0,2} =\mathbb Q$,
  so that $\mathcal K_2 = \mathcal K_1$. As a result
$|\Pi_0\cap \left(D_0\setminus D_1\right)| =
  |\Pi_1\cap \left(D_0\setminus D_1\right)| =\emptyset$.

  Further, $\widetilde H_{0,j} =L_j$ for $j\geq 3$ and
  $\widetilde H_{1,j} =L_j(\sqrt{2+t})$ for $j\geq 2$.
  Under the assumption of circular primitivity
  $[\widetilde H_{1,j}:\widetilde H_{0,j}]=[L_j(\sqrt{2+t}) : L_j] = 2$,
    for any $j\geq 2$. For $j=2$ it is self-evident and for $j\geq 3$
    the fields differ by their quadratic subfields,
    $\mathbb Q(\sqrt{2},\sqrt{2+t})$ and     $\mathbb Q(\sqrt{2})$,
    respectively. It follows that
    the sequence $|\Pi_j\cap \left(D_0\setminus D_1\right)|, j\geq 2$,
    is $2$-geometric.

    The partition inside $D_n\setminus D_{n+1}$ for $n\geq 1$
    is, up to scale, the same as in the generic case. Indeed, we can apply
    again Lemma 22 to conclude that $[\widetilde H_{n+1,j}:\widetilde H_{n,j}]= 2$
    for any $n\leq j-1, j\geq 3$.

    Finally taking into account that $|D_1| =1/2$ 
   we obtain 
   \[
   |\Pi_0(t)| =   |\Pi_1(t)| =0 +\frac{1}{2}\cdot \frac{1}{3}=
   \frac{1}{6},
   \]
   and the sequence $|\Pi_k(t)|, k\geq 2$, is $2$-geometric.
\end{proof}

To cover all cases we still need to address the failure of
circular primitivity. To that end let us observe that
for a circular $t$ and its associate value $w$ the respective
partitions are closely related. We have
$C_2(t)=t^2-2 = -(w^2-2)=-C_2(w)$, so that by Corollary 24
we get  
 \begin{equation}
 \begin{aligned}
 & \Pi_2(w)= \Pi_0(t)\cup\Pi_1(t),\ \ \ 
 \Pi_2(t) = \Pi_0(w)\cup\Pi_1(w),\\
 &\Pi_k(t) = \Pi_k(w), k = 3,4,\dots.
 \end{aligned}
 \end{equation}
 Let us note that the first two equations follow also directly from
 (6).

The phenomenon of elements of a partition being lumped together
and permuted as above can have arbitrary depth.
Indeed,  
let circular $t_0$ be circular primitive, $t_0^2+w_0^2=4$.
Let $2z_0=t_0+iw_0$, and for $N=2^k$ we consider
$2z_k=2z_0^N =C_{N}(t_0)+iw_0U_{N}(t_0) = t_k+iw_k$.
Clearly $t_k$ is circular
    and it is not  primitive, hence $w_k$, being also
    circular, is primitive and not circular primitive. 
    It follows from Corollary 24 and (7) that the partitions
    for $w_k$ and $t_0$ are connected as follows. 
\begin{equation}
 \begin{aligned}
 & \Pi_2(w_k)=\Pi_0(t_0)\cup\Pi_1(t_0)\cup\dots \cup\Pi_{k+1}(t_0),\ 
 \Pi_{k+2}(t_0) = \Pi_0(w_k)\cup\Pi_1(w_k),\\
 &\Pi_s(w_k) = \Pi_{k+s}(t_0), \ \ s = 3,4,\dots.
 \end{aligned}
 \end{equation}
 To get all the densities for $w_k$ it is sufficient to check
 that $|\Pi_0(w_k)|= |\Pi_1(w_k)|$. Assuming that we obtain 
  \[
  |\Pi_0(w_k)|= |\Pi_1(w_k)| =\frac{1}{3}\frac{1}{2^{k+1}},\ \
|\Pi_2(w_k)| =1-\frac{1}{3}\frac{1}{2^{k-1}}, 
  \]
  and the sequence $\{\Pi_s(w_k)\}_{s\geq 3}$ is $2$-geometric.
 Hence with increasing $k$ the set $\Pi_2(w_k)$ grows to include all primes,
while $|\Pi_0(w_k)\cup\Pi_1(w_k)|$ goes to zero.

 We failed to find a simple proof that $|\Pi_0(w_k)|= |\Pi_1(w_k)|$.
 We give here only a  sketch of the proof in the case $k=1$.
 It is based on the observation that the fields 
 $\widetilde H_{n,j}(w_1), n\leq j$ are real
 (which is always true in the circular case), and  
    if we adjoin $i$ to them we obtain the following radical
    extensions of $\mathbb Q(i)$. 
    Consider the binomials $z^{2^{j-2}}-i$ and $z^{2^n}-iz_0^2$. We claim that
    for $j\geq 2$ the field $\widetilde H_{n,j}(w_1)\left(i\right)$ is
    the splitting field of the
    two binomials over $\mathbb Q(i)$.

    Let us recall that in the circular case $\mathcal K_2 = \mathcal K_1$.
    It follows that there is no contribution into $\Pi_0$ and $\Pi_1$
    from $D_0\setminus D_1$ and $|D_1| = 1/2$, cf. Table 1. 
    The same holds for $D_1\setminus D_2$.
    To obtain that we need to inspect the fields $\widetilde H_{0,3}$
    and $\widetilde H_{1,2}$. These fields are the real parts of
    the splitting fields of $z^2 - i$,
    and of $z^2 = iz_0^2$, respectively. Clearly these
    splitting fields coincide, and hence their real part is equal
    to  $\widetilde H_{1,3}$. The consequence is that there is no
    contribution into $\Pi_0$ and $\Pi_1$ from $D_1\setminus D_2$
    and $|D_2| = 1/4$.

    We proceed to $D_2\setminus D_3$. Now we need to inspect the fields 
    $\widetilde H_{2,3}$, $\widetilde H_{2,4}$
    and $\widetilde H_{3,3}$. The fields are the real parts of
    the splitting fields of $z^2=i$ and $z^4 = iz_0^2$,
    of $z^4 - i$ and $z^4 = iz_0^2$,
    and of $z^2 - i$ and $z^8 = iz_0^2$, respectively.
    Clearly these three fields are all different so that
    $[\widetilde H_{3,4} : \widetilde H_{2,3}]=4$. Now we can claim that
    $|\Pi_0\cap\left(D_2\setminus D_3\right)| =
    |\Pi_1\cap\left(D_2\setminus D_3\right)|= |D_2|/4= 1/16$.

    It is tedious, but straightforward that the same claim
    will apply to $D_n\setminus D_{n+1}$ for all $n\geq 2$.
    Hence we obtain $|\Pi_0| = |\Pi_1|$.

    We believe that the same argument applies for any $k\geq 1$.
    However we found it too cumbersome to write down explicitly.

This completes the project of calculating the densities
of the partition for all values of $t$, except for the
reducible case which will be treated in the next section.

\section{Miscellanea}

A. The reducible case, $t^2-4 = a^2$ for a rational $a\neq 0$.

This case was studied first by Hasse, \cite{Hs},
then by Ballot,\cite{B1}. For completeness we  sketch here how it
can be treated along the lines of this work.
It is easier than the general case because now $\mathcal S_p$
is always of order $p-1$, (except when  $t^2-4 =0 \mod p$).

Section 2 applies without any changes. The proof of Theorem 5
contains everything necessary to obtain the new version in
the reducible case and it is actually simpler. 
\begin{theorem}
  For any prime $r\geq 2$ and a natural $j$
  the set $\mathcal K_j= \mathcal K_j(t,r)$
  contains exactly the primes
  $p\in \Pi(t)$  such that the polynomial $\Phi_j\mod p$ splits linearly.
\end{theorem}
The new version of Theorem 7 is also easier to prove.
\begin{theorem}
  For a prime $r\geq 2$ and any natural $n$ and $j\geq n$
  the set $\mathcal M_n\cap \mathcal K_{j}$
  contains exactly the primes
  $p\in \Pi(t)$  such that the polynomial $\Phi_{j}g_n \mod p$ splits linearly.
\end{theorem}
We denote by $\widehat H_{n,j}$ the splitting field
of the polynomial $\Phi_{j}g_n$.  In the reducible case $\xi$
is rational but the roots of the binomial $z^{r^n}-\xi$
are not under the assumption of $r$-primitivity. Let $z_n$ be a root
of this binomial. The $r$-genericity is void for $r\geq 3$.
We have the following version of Theorems 9 and 10 for $r\geq 3$.
  \begin{theorem}
    For any $r$-primitive $t$, $r\geq 3, t \neq 0,\pm 1, \pm 2$
    and any $n\geq 1$ and $j\geq n$
    \[
    H_{n,j}=\mathbb Q\left(\zeta_{j},z_{n}\right), \ \
    [\mathbb Q\left(\zeta_{j},z_{n}\right): \mathbb Q] = (r-1)r^{n+j-1}.
\]
  \end{theorem}
  By the methods used in the proof of Theorem 2 we can now easily
  establish its validity for $r\geq 3$ and a reducible $r$-primitive $t$.

 The reducible case for $r=2$ requires some attention,
  but it is still easier than the non-reducible case.
  The $2$-genericity is a legitimate condition and under this assumption
  Theorem 2 and the whole proof in Section 7
  is valid without any changes ($\sqrt{4-t^2}$ can be replaced
  with $i = \sqrt{-1}$). However, only
  non-genericity of type $A$ is possible. Again the result and the proof
  in this case are valid. Let us finally note that the condition 
  $2\pm t = 2a^2$ for some rational $a$ can be translated
  into the language of the rational $\xi$ where
  $t = \xi+\xi^{-1}$. We get that $t$ is non-generic of type $A$
  if and only if $\xi =\pm 2 b^2$ for some rational $b$.

  Let us note that the results about non-primitivity in
  Proposition 11 and in Corollary 24 apply in the reducible case
  without any changes.

  B, Sets of divisors for some sequences.

  Let us consider the sequences $\{w_n\}_{n\in \mathbb Z}$,
  $\{v_n\}_{n\in \mathbb Z}$ and $\{c_n\}_{n\in \mathbb Z}$,
  defined by the Chebyshev polynomials
    $w_n = W_{2n+1}(t), v_n= V_{2n+1}(t), c_n = C_n(t)$.
  We denote by 
  $W,V,C \in \mathcal \mathcal R(t)$ the matrices corresponding to
  respective sequences. We have $W=[-1\ 1], V= [1\ 1]$ and
  $C=[2\ t]$. Let us note that neither
  of these matrices belongs to $\mathcal S(t)$, because their
  determinants are $2+t, 2-t$ and $4-t^2$, respectively. \cite{L-W}.
  The sets of prime divisors of these
  sequences are equal to $\Pi_0(t),\Pi_1(t)$ and $\Pi(t)\setminus
  \left(\Pi_0(t)\cup\Pi_1(t)\right)$ for $r=2$.
  Indeed for any $k$ we have for the Lucas sequence $L_{2k+1} = Q^k W_{2k+1}$
  (\cite{L-W}, formula (1)). This clearly shows that the set of
  divisors for the sequence $W$ is equal to $\Pi_0(t)$ for $r=2$.
  To get the second claim  we can now employ the twin symmetry
  $\Pi_1(t) = \Pi_0(-t)$ and the identity $V_{2k+1}(t) = U_{k+1}(t)-U_k(t)
  =(-1)^k W_{2k+1}(-t)$.

  For the last claim we use the formula $L_{4k}= c U_{2k}(t)
  = cU_{k}(t)C_k(t), c = TQ^{2k-1}$.

  In the cubic case $t^2-4=-3b^2$ we consider the sequence
  $S = \frac{1}{2b}[2\ t+b] $. 
  The set of divisors of $S$ is equal to
  $\Pi(t)\setminus \Pi_0(t)$ for $r=3$.
  We provide
  a proof in the language of the present paper. We use the fact that
  $S$ has determinant equal to $1$, and so it is an element of $\mathcal S(t)$
  and $S^3=I$.   The prime $p$ belongs to  $\Pi(t)\setminus \Pi_0(t)$
  if and only if $3 \mid \chi(t,p)$. It follows that for $k =\chi(t,p)/3$
  we have $D^{k} = S\mod p$ or  $D^{k} = S^{-1} \mod p$. Accordingly
  $SD^{-k} = I \mod p$ or  $SD^{k}=I \mod p$, which shows that.  
  there is a zero in the sequence $S\mod p$.

  In the cubic primitive
  case the formulas from Theorem 13 apply and give us the
  prime density for the divisors of $S$ equal to $3/4$.
  It was proven in \cite{L} in a special case,
  and in \cite{B3} in general under genericity
  assumptions that seem to be different from ours.

  We turn our attention to the sequences $WS, VS$ and $CS$.
  The new feature is that they are not elements of $\mathcal S(t)$.
   It was established in \cite{L-W}, Proposition 30,
   that 
  their sets of divisors contain primes $p$ such that 
  $\chi(t,p) = 3\mod 6$, 
  $ \chi(t,p) = 6 \mod 12$ and $ 12\mid \chi(t,p)$, respectively.
  In particular these sets are
  disjoint and their union is equal to $\Pi(t)\setminus \Pi_0(t)$
  for $r=3$.

  In \cite{B3} Ballot calculated that, under some genericity assumptions,
  the density of primes for which the index is divisible by $6$,
  is equal
  to $1/2$. Using our method we can prove a little bit more.
  \begin{theorem}
    For a cubic primitive and twin primitive $t$ the sets of divisors
    of the sequences $WS, VS$ and $CS$ have densities equal to $1/4$.
\end{theorem}
  In view of the above characterisation of the sets of divisors
  the theorem is equivalent to the following statement.
  \begin{theorem}
    For a cubic primitive and twin primitive  $t$
    the following sets of primes have densities 
\[
|\left(\Pi(t)\setminus\Pi_0(t,3)\right) \cap \Pi_0(t,2)| =
|\left(\Pi(t)\setminus\Pi_0(t,3)\right) \cap \Pi_1(t,2)|=\frac{1}{4}.
\]
    \end{theorem}
  \begin{proof}
    We first claim is that $t$ is automatically $2$-generic.
    Clearly $t$ cannot be circular or non-generic of type $B$.
    Suppose that our $t$ is non-generic of type $A$, i.e.,
    $2-t = 2e^2$ for a rational $e$. It follows that $2+t=
    (4-t^2)/(2-t)=6f^2$ for a rational
    $f$. To see that these equalities cannot be satisfied simultaneously
    we employ the representation of $t$ by two integers
    $T,Q\neq 0$. We have $2+t =T^2/Q$ and hence $Q=6F^2$ for a natural $F$.
    Further $2e^2=2-t = 4-T^2/Q$ so that $T^2= 6(4-2e^2)F^2$.
    The last equality can be turned into the diophantine equation
    $2A^2=3B^2+C^2$ for integer $A,B,C$. This equation has no integer
    solutions. Indeed if $B$ and $C$ are even then $A$ must be even
    as well, which allows cancelation of the factor $4$ on both sides.
    After further cancelations, if needed,
    we end up with the same equation with
    both $B$ and $C$ odd. In such a case $3B^2+C^2$ is divisible by
    $4$ and not divisible by $8$, which cannot hold for $2A^2$.

    The second step involves a new characterization of the set
    $\mathcal K_j\cap M_n$ for $r=3$, namely it is the set of primes $p$ such
    the polynomial
    $ h_{j,n}(x)= \left(C_{3^{j-1}}(x)+1\right)
    \left(C_{3^n}(x)-t\right) \mod p$
    splits linearly. It is based on the observation that the special
    sequence group $\mathcal S(t)$ contains the third root of identity
    $S$ with the trace equal to $-1$. It follows by Lemma 8 that
    if $\left(C_{3^{j-1}}(x)+1\right)  \mod p$ has a root $\eta\in \mathbb F_p$
    then $\mathcal S_p(t)$ contains $Z_j$ with the trace $\eta$ which is
    by necessity a root of identity of degree 
    $3^j$. The rest of the argument is the same as in the proof of
    Theorem 7.
    By the proof of Theorem 13 (which holds also for $r=3$
    under the assumption of cubic primitivity)
    we know the density of this set, namely $1/3^{n+j-1}$.
    We conclude that the Galois group of the splitting field of
    $h_{j,n}$ is of order $3^{n+j-1}$.

    By Theorem 18 we have that $\mathcal K_s\cap M_m$ for $r=2$
    is characterized as the set of primes $p$ such that the polynomial
    $C_{2^{s-2}}(x)\widetilde f(x) g_m(x)$ splits linearly. By Theorem 21 the 
    Galois group of the splitting field has the order $2^{m+s-1}$.
    Since this order has no common divisors with the order of the previous
    Galois group we conclude that the Galois group of the joint
    splitting field of the polynomial $h_{j,n}C_{2^{s-2}}\widetilde f g_m$
    is the direct sum of the two Galois groups, with order equal to
        $3^{n+j-1}2^{m+s-1}$.

    For any $j,n, s, m$ the intersection of the sets 
        $\mathcal K_j\cap M_n$ for $r=3$ and 
    $\mathcal K_s\cap M_m$ for $r=2$ is characterized by linear splitting
    of the polynomial $h_{j,n}C_{2^{s-2}}\widetilde f g_m$.
    Further the density of the intersection is the product of the densities
    of the sets. By the elementary probability theory it follows that
    the divisibility properties  by the powers of $3$ and powers of  $2$
    of the index $\chi(t,p)$
    are
    ``independent'', in particular 
    \[|\left(\Pi(t)\setminus\Pi_0(t,3)\right) \cap \Pi_0(t,2)| =
     |\Pi(t)\setminus\Pi_0(t,3)| |\Pi_0(t,2)| =
     \frac{3}{4}\frac{1}{3} = \frac{1}{4}.
     \]
        \end{proof}
  Let us note that it follows from the proof of Theorem 7
  that in the cubic case the splitting field of the DeMoivre's polynomial
  $C_{3^n}(x)-t$ is equal to the splitting field of the polynomial
  $h_{n,n}$. Girstmair, \cite{G}, studies the Galois groups of
  DeMoivre's polynomials for odd degrees. The paper contains 
  the claim, without all the details of the proof,
  that for a cubic primitive $t$ the order of the
  Galois group of the splitting field of $C_{3^n}(x)-t$
  is equal to $3^{2n-1}$.


  C. Linear second order recursive sequences with the set of divisors
  equal to $\Pi_0(t,r)$ for any $r\geq 3$ and
  $\Pi_0(t)\cup\Pi_1(t)$ for $r=2$..

  For a given $t$ and a prime $r$  let $t_r = C_r(t)$.
  By Lemma 8 there is a matrix $Y_1 \in \mathcal S(t_r)$ such that $tr\ Y_1 =t$.
  The sequences corresponding to $Y_1$, and also the powers
  $Y_m = Y_1^m, m=2,\dots, r-1$, have the same sets of prime
  divisors equal to $\Pi_0(t,r)$. To prove it we invoke the isomorphism
  of $\mathcal R_{t}$ and $\mathcal R_{t_r}$ from \cite{L-W}, used already in the proof
  of Proposition 11. This time we need to recall the interpretation
  of the isomorphism $\Psi_r: \mathcal R_{t} \to R_{t_r}$ in terms of
  sequences rather than matrices. Namely it was observed in \cite{L-W}
  that for a recursive sequence $\{ x_n \}_{n\in\mathbb Z}$ for the parameter
  $t$ the sequences  $\{ x_{kr+s} \}_{k\in\mathbb Z}, s = 0,\dots r-1$,
  are recursive sequences for the parameter $t_r$. Further
  $\Psi_r([x_0\ x_1]) = [x_0\ x_r]$ and $\Psi (D_t^r) = D_{t_r}$.

  We choose $Y_1 =\Psi (D_t)$, and hence $Y_s= \Psi (D_t^s)$. In terms of
  sequences, if $u_n = U_n(t), n\in \mathbb Z$, is the identity sequence
  in $\mathcal R(t)$ then clearly the sequences
  $\{ u_{kr+s} \}_{k\in\mathbb Z}, s = 0,\dots r-1$, represent the elements
  $Y_s \in \mathcal R(t_r)$. Let us restrict our attention to the sequence
  $\{ u_{kr+1} \}_{k\in\mathbb Z}$. Using  the formula for the Lucas sequence
  $L_{2m} = c U_m(t), c = TQ^{m-1}$,
    we obtain by (1) that
    if $ u_{kr+1}=0 \mod p$ then 
    $\chi(p,t)|2(kr+1)$, and hence $p\in \Pi_0(t)$ for $r\geq 3$, and
    $p\in \Pi_0(t)\cup \Pi_1(t)$ for $r=2$.

    Conversely, for $r\geq 3$, if 
  $p\in \Pi_0(t)$ then $\chi(t,p) = mr+s$ for some integer  $m$ and
  $s, 1\leq  s\leq r-1$. Further there is a natural $l$ such that
    $ls=qr+1$ for a natural $q$, which delivers
    $u_{l(mr+s)} = u_{kr+1} =0 \mod p$
    for $k=lm+q$,
  so that indeed $p$ is the divisor of the sequence $Y_1$.
  It is also clear from the proof
  that the sequences $Y_1^s, s = 0,\dots r-1$ have the same
  set of prime divisors, independent of $s$.

     For $r =2$, if 
     $p\in \Pi_0(t)$ then $\chi(t,p) = 2m+1$ for some integer $m$,
     and hence $u_{2m+1} =U_{\chi}(t) = 0 \mod p$. If
     $p\in \Pi_1(t)$ then $\chi(t,p) = 2(2m+1)$ for some integer $m$,
     and hence $U_{2(2m+1)}(t) = U_{2m+1}(t)C_{2m+1}(t)= 0 \mod p$.
     The divisors of $C_{2m+1}(t)$ belong to $\Pi_2(t)$, hence
     $u_{2m+1} = U_{2m+1}(t) =0 \mod p$. In both cases $p$ divides the
     sequence $Y_1$.

  D. Ballot's sequence   $B_k= L_{rk}/L_k$.

  For the convenience of the reader we provide here the proof 
  that the set of prime divisors of  the sequence
  $B_k= L_{rk}/L_k$ is equal to $ \Pi(t) \setminus \Pi_0(t,r)$.

  Indeed, if   $\chi(t,p) = rn$ then $p \mid L_{rn}$ and
  $p\nmid L_n$, so that clearly $p\mid B_n$.

  To prove the opposite we observe that 
  if a prime $p$ divides the numerator and does not divide the denominator
  then   $\chi(t,p)$ is divisible by $r$. To finish the proof
  we will check that, since by assumption $p\neq r$, 
  if $p$ divides the denominator $L_k$  then it does not divide
  $B_k$.  
  To that end we use the expression of the Lucas sequence by the
  Chebyshev polynomials, formulated in the Introduction.
  In the case of prime $r\geq 3$ we get for 
  $k$ odd and for $k=2s$, respectively, 
  \begin{equation}
  B_k= Q^\alpha \frac{W_{rk}(t)}{W_k(t)} = Q^\alpha W_r(C_k(t)), \ \
  B_k= Q^\alpha \frac{U_{rs}(t)}{U_{s}(t)} = Q^\alpha U_r(C_s(t)),
  \end{equation}
  where  $\alpha = (r-1)k/2$.
  These equalities are of independent interest and they
  follow from the identities
  $U_{mn}(t) = U_m(C_n(t))U_n(t)$ and
  $W_{mn}(t) = W_m(C_n(t))W_n(t)$, for any $m,n$, and odd $m,n$,
  respectively, cf. \cite{W}. 
  In particular, since both $W_r(C_k)$ and $U_r(C_s)$ have degrees
  equal to $(r-1)k/2$, we can see that $B_k$ are integer.

  If $p$ divides $W_k(t)$ then by the identity
  $(t-2)W_k^2(t) = C_k(t) - 2$, following from (1),
  we conclude that $C_k(t) = 2 \mod p$. It follows that
  $W_r\left(C_k(t)\right) = W_r(2) \mod p$.
  It is easy to check that $W_r(2) = r$,
  \cite{W}. Similarly  if $p$ divides $U_s(t)$ then by (2) 
  $ C_s(t) = \pm 2 \mod p$. Since  $U_{n}(\pm 2) = \pm n$
  we conclude that $U_r(C_s(t)) = \pm  r\mod p$.
  Hence in both cases if $p$  divides the denominator $L_k$
  then it does not divide $B_k$,
  which ends the proof for $r\geq 3$.

  There is a version of (9) in case of $r=2$. Indeed we get
  for $k$ odd and for $k=2s$, respectively, 
  \begin{equation}
  B_k= Q^\alpha \frac{U_{k}(t)}{W_k(t)} = Q^\alpha V_k(t), \ \
  B_k= Q^\beta \frac{U_{2s}(t)}{U_{s}(t)} = Q^\beta C_s(t),
  \end{equation}
  where  $\alpha = (k-1)/2$ and     $\beta = s$.
  We can see again that the sequence $B_k$ is integer,
  and that its set of prime divisors is equal to the union
  of such sets for the sequences $V$ and $C$ discussed in
  part B of this Section. Our claim follows immediately.

  \section{Applications to arithmetic dynamics:
    the sets of divisors of iterates of linear and quadratic functions}

A. Rotations.

For a rational $t,  -2 <t < 2, t\neq 0,\pm 1$, the linear recursive sequences:
$x_{n+1}= t x_n-x_{n-1}, \ n\in \mathbb Z$,
 can be interpreted as orbits of an irrational rotation 
 \[
   [x_{n+1}\ x_{n+1}] =[x_0\ x_1]\left[ \begin{array}{cc} 0 &     -1\\
    1  &   t \end{array} \right]^n = [x_0\ x_1]D^n,
\]
on invariant ellipses $x_1^2-tx_1x_0+x_0^2 = c$, for a rational $c$.
It was discovered in \cite{L-W} that for initial conditions $[x_0\ x_1]$
with $c$ which is not a rational square
the set of divisors of the recursive sequence
is not bigger than half of all prime numbers.
Indeed, if $x_n = 0 \mod p$ then 
\[
c = x_{n+1}^2-tx_{n+1}x_{n}+x_n^2= x_{n+1}^2 \mod  p.
\]
It is well known that if $c$ is not a rational square
then the set of primes such that
$c \mod p$ is a square has prime density $1/2$.

The numerical results in \cite{L-W} indicate
that the prime density for generic initial conditions 
is between $0.32$ and $0.36$.
Exact results are summarized 
in the survey paper
of Moree, \cite{M3},Chapter 8.4.

It was proven  that with the exception of $t=0,\pm 1, \pm 2$
the set of divisors is infinite for any rational non-zero
initial conditions, \cite{Polya}, \cite{Wr}, \cite{Lw}.
Schinzel, \cite{S2}, proved in the reducible case
that there are also infinitely many non-divisors.
With our method we can strengthen this result to
general case and positive density of non-divisors.
It is presented in Section 12.

B. Chaotic dynamics, iterations of $C_k(x)$.

In the paper \cite{J} Rafe Jones established that
the set of prime divisors of orbits
for  several families of quadratic polynomials
for integer initial conditions has density zero. The Chebyshev  
polynomial $C_2(x) = x^2-2$ is a special case in his theorem. 
This result can be extended to any Chebyshev polynomial $C_k(x), k\geq 2$,
and all rational initial conditions $x_0 \neq 0,\pm 1,\pm 2$.
The resulting set of divisors for the numerators of the orbit
has density zero. Indeed, the $n$-th composition of
$C_k(x)$ is equal to $C_{k^n}(x)$. For simplicity let us choose
a prime divisor $r\geq 2$ of $k$ and let $r^s||k$. It is immediate that if
$C_{k^n}(x_0) = 0 \mod p$ for some natural $n$
then $C_{4k^n}(x_0) = 2 \mod p$, and hence
$\chi(x_0,p) = 4k^n$ and $p \in \Pi_{ns}(x_0,r)$, or
$p \in \Pi_{ns+2}(x_0,r)$ for $r=2$.
Incidentally, it follows that the
numerators of the orbit are relatively prime (except  for the
possible common factors
$2$ or $r$). Moreover, since $\chi(x_0,p)| p \pm 1$ we get that
any prime divisor $p$ of the $n$-th element of the orbit
is bigger or equal than  $4k^n-1$. To finish the proof of density $0$
let us note that the set  $\mathcal A$ of all the prime  divisors of the
numerators can be split as follows
    \[
\mathcal A =     \left(\mathcal A\cap\bigcup_{n=0}^N\Pi_{ns}(x_0,r)\right)
\cup
    \left(\mathcal A\cap\bigcup_{n=N+1}^{\infty}\Pi_{ns}(x_0,r)\right).
    \]
    The first set is finite and the second, being
    a subset of a set with a small prime density by Theorem 2,
    has small upper prime density. Increasing $N$ we obtain zero prime
    density for $\mathcal A$.

    C. Other quadratic iterations.

    In the paper
``Quadratic Recurrences with a Positive Density of Prime Divisors''
Gottesman and Tang, \cite{G-T}, study the iterations of
the quadratic function
$f_t(x) = (x+t)^2-2-t$ at $x=0$, for an integer $t$.
They also point to the conjectural connection with
the iterations of $g_s(x) = (x+s)^2-s$. Our results
shed some light on this issue.

By the change of variable
$y=x+t$ we obtain that $f_t$ is conjugate to the
Chebyshev polynomial $C_2(y) = y^2-2$. Similarly  $g_s$ is conjugate
to $y^2$ by the change of variable $y=x+s$. Hence for the iterations 
we get $f_t^n(x) = C_{2^n}(x+t)-t$ and $g_s^n(x) = (x+s)^{2^n}-s$.
For a rational $t$ the set of divisors  of the sequence of
values of $f_t^n(x)$ and $g_s^n(x)$ at $x=0$ are  equal to the set of
divisors of
the sequences  $C_{2^n}(t)-t$ and $s^{2^n}-s$, respectively.
The second case can be cast in the same form as the first.
Indeed, if we let $t = s+s^{-1}$ then we get that
the sets of prime divisors of the sequences
$s^{2^n}-s$ and $C_{2^n}(t)-t$ coincide, with the exception of the divisors
of the numerator and denominator of $s$. Hence the problem for
$g_s$ becomes a special case of the problem for $f_t$,
namely the case of reducible $t$.

For a given prime
$p$ let us invoke the special sequence group $\mathcal S_p(t)$.
It follows from  Lemma 8 that $C_{2^n}(t)-t =0 \mod p$ 
if and only if $D_t^{2^n} =D_t \mod p $ or $D_t^{2^n} =D_t^{-1}\mod p$. Further,
it is equivalent to $D_t^{2^n\pm 1} =I \mod p$. Since every
odd prime divides some element of the sequence $2^n-1$,
the last property is equivalent to the index $\chi(t,p)$ being odd.
Hence the set of divisors of the iterations of $f_t(x)$ or $g_s(x)$
at $x=0$ is  equal to
$\Pi_0(t)$ for $r=2$.

In this way our results deliver the content of Theorem
1.1 in \cite{G-T}, with the improvement that the condition
$t^2-4 \neq 2b^2$ for a rational $b$ is redundant. The condition for
non-genericity
of type $B$ says    $t^2-4 = -2b^2$ for a rational $b$,
which is void  in the case of integer $t$
studied in \cite{G-T}.

For the iterations of $g_s(x)$ at $x=0$ we get
only the reducible case and as explained in the previous section
the condition of twin primitivity is that $s\neq \pm b^2$ for a rational
$b$. Further, the condition of non-genericity of type $B$ is that
$s = \pm 2b^2$ for a rational $b$.

\section{The sets of non-divisors}

We have shown in Section 11, part A, that for any 
nonzero sequence $X\in R(t)$ with determinant which is not a rational square
the set of non-divisors
contains the set of primes $p$ such that $\det X$ is
not a square $\mod p$. This set has density $1/2$,
cf. \cite{L-W}, Theorem 31.
If $\det X = c^2$ for a rational $c$ then $Y=c^{-1} X$ belongs
to $\mathcal S(t)$. We will prove in this Section
that also for any such sequence, with the exception of sequences
with a zero element,
there is a  set of non-divisors with positive density.
A sequence  $Y\in \mathcal S(t)$ has a zero element if and only if
$Y= \pm D^k$ for some integer $k$.

Let us first inspect the case when $Y=AD^k$ where $A$ is an element of
finite order in $\mathcal S(t)$ and $k\in \mathbb Z$. It is well known that
except  for $\pm 1$ there are elements of finite order in a quadratic field
only for the cyclotomic quadratic fields $\mathbb Q(i)$ and
$\mathbb Q(\sqrt{-3}$. In our terminology this is the case of circular 
or cubic $t$.

In the circular case we need to consider the
set of prime divisors of $C\in \mathcal S(t)$ such that $tr\ C =0$.
For circular $t$ there is essentially only one such element
and In Section 10, part B,
we established that its set of non-divisors is equal to
$\Pi_0(t)\cup \Pi_1(t)$ for $r=2$. In the circular case
the density of this set is by Theorem 28 equal to $1/3$ if $t$
is circular primitive. If $t$ is not $2$-primitive then density
of $\Pi_0(t)\cup \Pi_1(t)$ exceeds $1/3$ by Proposition 11.
If $t$ is $2$-primitive but not circular primitive then by the formulas
(8) this density
is equal to $\frac{1}{3}\frac{1}{2^k}$ for some $k\geq 1$. Note that this
claim  does not depend on the equality $|\Pi_0(t)|=|\Pi_1(t)|$,
which has a long and complicated proof.

  In the cubic case we need to consider the set of prime divisors of
  the element $S$ of order $3$ in $\mathcal S(t)$. By the results 
  of Section 10, part B, the set of non-divisors
  is equal to $\Pi_0(t)$ for $r=3$. Its density in the cubic  primitive case
  is equal to $1/4$ by the formulas in Theorem 13. If $t$ is not $3$-primitive
  then by Proposition 11 the density of the  set of non-divisors exceeds
  $1/4$. If $t$ is $3$-primitive but not circular primitive then
  it follows from the formulas (5) that the density of the set of non-divisors
  is equal to $\frac{1}{4}\frac{1}{3^k}$ for some $k\geq 1$.

  We proceed to the general case with the new assumption that
  $Y \neq A D^k$ for any torsion element $A \in \mathcal S(t)$
  and any $k\in\mathbb Z$. In particular $b= tr\ Y \neq 0,\pm 1, \pm 2$. 
For such a sequence $Y$ we consider the  order 
$ord_p(Y)$ of $Y\mod p$ in  $\mathcal S_p(t)$. To deal with two sequences
$D$ and $Y$ simultaneously we remove from considerations
the prime divisors of the denominators of  $t, b$ and all the
elements of $Y$. 
It is straightforward that
$ord_p(Y) = \chi(b,p)$. Indeed, by Lemma
 7 the rings $\mathcal R(t)$ and $\mathcal R(b)$ are
isomorphic, which leads to the isomorphisms of
the respective sequence groups. 
Moreover
an odd prime $p$ divides the sequence
$Y\in \mathcal S(t)$ if and only if $ord_p(Y)$ divides
$2\chi= 2\chi(t,p) = 2  ord_p(D)$.
Indeed, the divisibility is equivalent to the existence
of an integer $k$ such that $D^kY = \pm I \mod p$. 
The claim now follows from the following fact:
$ord(h)\mid ord(g)$ for two elements $h,g$ of a finite cyclic group
if and only if $h$ belongs to the subgroup generated by $g$.

We will need the following lemma
  which will be proven at the end. 
  \begin{lemma}
    For $t\neq 0,\pm 1,\pm 2$, if a sequence $Y\neq AD^k$
    for any torsion element $A\in\mathcal S(t)$ and any integer
    $k$, then for all but finitely many
    primes $r$ both $t$ and $b = tr\ Y$ are $r$-generic
    and $YD^k \neq Z^r$ for any integer $k$ and $Z\in \mathcal S(t)$. 
     \end{lemma}
\begin{theorem}
  For any $t\neq 0,\pm 1,\pm 2$ 
  and a sequence  $Y\in \mathcal S(t)$ with no zero elements,
  there is a subset  of positive density of primes which do not divide
  any element of $Y$.
\end{theorem}
\begin{proof}
  Thanks to the completed discussion of $Y$ equal to torsion elements 
  we can assume that the assumptions of Lemma 34 are satisfied
  and we fix the smallest prime $r\geq 5$ with the properties
  described there. .

  Clearly the set of primes $p$ such that $r\mid ord_p(Y)$ and
  $r\nmid \chi(t,p)=ord_p(D)$
  is contained in the set of non-divisors of $Y$.
  We further consider the set
  $B_1 = \left(\mathcal K_{1}\cap\mathcal M_1\right)\setminus
  \mathcal K_{2}$ from the first table in Table 1, which contains only
  primes with the index non-divisible by $r$.

  We can repeat the constructions of Section 2 replacing $D$ with $Y$,
  and obtain similar conclusions.
  Thus we introduce $\mathcal M_1(b) =
  \{ p \in \Pi(t)\ | \ \exists G\in \mathcal S_p(t),  G^{r}=Y\}$.
  We have that $E(b)=\mathcal K_{1}\setminus\mathcal M_1(b)$ contains 
  only primes such that $ord_p(Y)$ is divisible by $r$.

  While the proof of Theorem 2 delivers the densities
  $|B_1| = \frac{1}{r^2} , |E(b)| = \frac{1}{r}$, 
  it remains to calculate the density of the intersection.
  To that end let us observe that we need to consider only the
  four subsets
  \[
    Z_j=\mathcal K_{j}\cap\mathcal M_1, \ \
    \widehat Z_j =\mathcal K_{j}\cap\mathcal M_1\cap\mathcal M_1(b), \ j=1,2.
    \]
    Indeed we have that $|B_1\cap E(b)| =
    \left( |Z_1| - |\widehat Z_1|\right) - 
    \left( |Z_2| - |\widehat Z_2|\right)$.
    Since we know the densities of the sets
    in Table 1, it follows that $|Z_1| = \frac{1}{(r-1)r},
    \ |Z_2| = \frac{1}{(r-1)r^2}$.

    We proceed with the characterization of the sets $\widehat Z_j, j=1,2$,
    by the methods of Sections 3 and 4.
    We obtain that an odd  prime $p$ belongs to
    $\widehat Z_j, j=1,2$, if and only if the polynomial $f\Phi_j\mod p$
    splits linearly or splits quadratically and  the DeMoivre's
    polynomials
    $C_r(x)-t\mod p$ and $C_r(x)-b\mod p$ split linearly.
    By Theorem 9  the splitting field $H_{0,j}$ of the polynomial     $f\Phi_j$
    is equal to $\mathbb Q(\xi, \zeta_j)$, where $\xi$ is a zero of $f$ and
    $\zeta_j$ is a primitive root of $1$ of order $r^j$. 
    The splitting field of the polynomial 
    $f\Phi_j\left(C_r(x)-t\right)\left(C_r(x)-b\right)$
    is the radical extension
    $H_{0,j}\left( \xi^{\frac{1}{r}},\beta^{\frac{1}{r}}\right)$,
          where $\beta$ is equal to $Y$ under the identification
          of $\mathcal R(t)$
          with $\mathbb Q(\xi)$, (recall that with this identification $\xi$
          is equal to $D$).
        By our assumption that both $D$ and $Y$ do not have an $r$-th root
        in $\mathcal R(t)$  we get that the binomials
        $x^r - \xi$ and $x^r - \beta$
        are irreducible over $\mathbb Q(\xi)$
        and hence also over its abelian extensions 
        $\mathbb Q(\xi, \zeta_j), j= 1,2$. In the last claim
        we take advantage of the fact that by the choice of $r$ the field
        $\mathbb Q(\xi)$ does not contain an $r$-th root of unity.

        It remains to show that the intersection of the extensions
        $H_{0,j}\left( \xi^{\frac{1}{r}}\right)$
          and         $H_{0,j}\left( \beta^{\frac{1}{r}}\right)$
            is equal to the field $H_{0,j}$. By the theorem of Schinzel,
            \cite{S1}, the two fields coincide if and only if
            $\beta = \xi^k\alpha^r$ for an integer $k, r\nmid k$,
            and $\alpha \in H_{0,j}$.
            By this condition             $\alpha^r=\beta\xi^{-k}$
            is an element of $\mathbb Q(\xi)$ and we can argue that
            also $\alpha\in \mathbb Q(\xi)$.
            Indeed, if $\alpha\notin \mathbb Q(\xi)$ then the binomial
            $x^r-\alpha^r$ is irreducible over $\mathbb Q(\xi)$,
            and hence also over its abelian extension $H_{0,j}$.
            Since $\alpha \in H_{0,j}$ we obtain a contradiction.

            Translating the condition  $\beta = \xi^k\alpha^r$
            into $\mathcal R(t)$ we obtain 
            $Y =D^kA^r$, where $A$ corresponds to $\alpha$ under
            the isomorphism of $\mathcal R(t)$ and $\mathbb Q(\xi)$.
            This is excluded by the choice of $r$.
            The conclusion is that 
            $[H_{0,j}\left( \xi^{\frac{1}{r}},
              \beta^{\frac{1}{r}}\right):H_{0,j}] = r^2$,

            Now the proof of Theorem 10 can be used to obtain the
            structures and the orders of the respective Galois groups.
            By the arguments employed in Section 4 we obtain then the densities
            $|\widehat Z_1| = \frac{1}{(r-1)r^2},
            \ |\widehat Z_2| = \frac{1}{(r-1)r^3}$,
            which gives us immediately $|B_1\cap E(b)| =
            \frac{r-1}{r^3}$.
            \end{proof}
In giving the exact density of the subset $B_1\cap M_1(b)$ the
proof delivers  more than needed, but we can use it in the study of the
sets of divisors of the powers $Y^{2^n}$. Assuming that $t$ is
$3$-generic (i.e., $3$-primitive and not cubic), and $YD^k$ is not
a third power for any integer $k$, we get that $r=3$ can be used in the proof
with $Y$ replaced by any power $Y^{2^n}$. Hence all these sequences
for any natural $n$ have a set of non-divisors of density $2/27$.
Note that the sets are actually  the same for different powers $2^n$ because
the condition that $Y\mod p$ has a third root in $\mathcal S_p(t)$
is equivalent to $Y^{2^n}\mod p$ having a  third root there. Hence the union
of the sets of divisors of all the powers 
$Y^{2^n}$ is disjoint from the set of density $2/27$.

We can repeat this argument
for the sequence of powers $Y^{k^n}$ for any fixed natural $k$,
with an appropriate choice of a prime $r\nmid k$.
We get a common subset of non-divisors of positive density.
It remains an open problem if the densities (or at least upper densities)
of the sets of prime divisors
of the powers $Y^n$   converge to $1$ as $n\to\infty$.

The last theorem holds also in the reducible case
by the modifications described in Section 10.A.
It is then a strengthening of the theorem of Schinzel,
\cite{S2}, who proved that for sequences
$\{a^n-b\}_{n\in \mathbb N}$, with integer $a,b$ and $b\neq a^k$ for any
$k$, the set of non-divisors is always infinite.

We conclude with the proof of Lemma 34.
\begin{proof}
  We consider the equation $YD^k = Z^r$ for a prime $r\geq 5$
  where $Z$ is an unknown element of the group $\mathcal S(t)$
  and $k$ is an unknown integer. We will be proving that
  this equation has no solutions for all but finitely many  prime $r$.

  Let $\mathcal T(t)$ be the torsion subgroup of $\mathcal S(t)$
  and 
  $F: \mathcal S(t) \to \widetilde{\mathcal S(t)}
  = \mathcal S(t)/\mathcal T(t)$ be the canonical group homomorphism.
  Let $\widetilde D =F(D), \widetilde Y =F(Y)$
  and $\widetilde Z =F(Z)$.
  
  The torsion elements in $\mathcal S(t)$ are in general only $\pm I$.
  In the circular and cubic case there are some other torsion elements
  in $\mathcal S(t)$ of order $2$ and $3$ only.
It follows that for
  $r\geq 5$ the torsion elements in $S(t)$ are trivially 
 $r$-th powers. Hence the equation $YD^k = Z^r$ is
solvable in $\mathcal S(t)$ if and only if the
equation $\widetilde Y{\widetilde D}^k = {\widetilde Z}^r$ is solvable 
 in the quotient group $\widetilde{\mathcal S(t)}$.
 Since $ \mathcal S(t)$ is a subgroup of the multiplicative
 group of a quadratic number field we can conclude that
 $\widetilde{\mathcal S(t)}$ is a countable free abelian group. 
 Both $\widetilde D$ and $\widetilde Y$ are non zero elements in this
 group
 by our assumptions. By the properties of free abelian groups
 we have that there is exactly one natural number
 $m$ such that $\widetilde D = {\widetilde R}^m$, for some $\widetilde R
 \in \widetilde {\mathcal S(t)}$ which has no roots. We obtain an
 equivalent equation :  $\widetilde Y\widetilde {R}^{km} = \widetilde {Z}^r$.

 Since  
 $\widetilde R$ is not any power it can be chosen as the first
 element of a basis in $\widetilde{\mathcal S(t)}$.
 We consider the quotient group
 $\widehat{\mathcal S(t)} =\widetilde{\mathcal S(t)}/<R>$,
 where $<R>$ is the cyclic subgroup generated by $R$.
 It is clearly also a countable free
 abelian group. If our equation has a solution in
  $\widetilde{\mathcal S(t)}$ then also the
 equation $\widehat Y={\widehat Z}^r$
 has a solution in $\widehat{\mathcal S(t)}$,
 where $\widehat Y$ is the homomorphic image of
 $\widetilde Y$ in the factor group and
 $\widehat Z$ is the unknown element of the
 factor group.

 In a countable free abelian group any element different from
 the unit element can be an $r$-th power
 for only finitely many primes $r$. It remains to show that
 $\widehat Y$ is not the unit element. Should that be the case
 we would have $D= A_1R^m, Y=A_2R^{km}$ for some $R\in \mathcal S(t)$,
 some torsion elements
 in $A_1, A_2\in \mathcal S(t)$ and integer $m,k$.
 These relations lead to $Y = AD^k$ for some torsion element $A$,
 the case that was treated at the beginning of this Section.

 It remains to check that for almost all $r\geq 5$ both
 $t$ and $b$ are $r$-generic. We need only to guarantee that
 both $t$ and $b$ are $r$-primitive for almost all $r$.
 The remaining condition can be violated for only two values of $r$,
 one for $t$ and one for $b$. 
 We have established that
 $D = AR^m$, where $A$ is a torsion element and $R$ has no roots,
 and now if $t$ is not $r$-primitive then  also $D = P^r$, for some $P$.
 Since we can absorb $A$ into $P$ we get 
 $P^r = R^m$. Taking this equality homomorphically into the
 factor group $\widetilde {S(t)}$ we would obtain that $r\mid m$.
 There are only finitely many divisors a fixed natural number $m$.
 
 If $b$ is not $r$-primitive then $Y= Z^r$, for some $Z$,
 which we have already established
 can happen for only finitely many $r$. 
\end{proof}

{\bf Acknowledgments}
The author is grateful to Andrzej D\c{a}browski and Zbigniew Lipiński
for many fruitful discussions. We would also like to acknowledge
the hospitality of the Center for Geometry and Physics
at Stony Brook, where part of this paper was written during
the visit in November-December of 2023.
The final version of this paper
was influenced by our friendly discussions with
Dale Brownawell and John Lesieutre
at the Algebra and Number Theory Seminar at the Department of
Mathematics at Penn State
University in December of 2024.

\

\end{document}